\magnification=1150
\overfullrule=0pt
\parindent=0pt

\input amssym.def
\input amssym.tex
\input epsf

\def\pec{\vskip2mmplus1mmminus1mm}
\def\gec{\vskip4mmplus1mmminus1mm}

\def\og{\leavevmode\raise.3ex\hbox{$\scriptscriptstyle\langle\!\langle\,$}}
\def\fg{\leavevmode\raise.3ex\hbox{$\scriptscriptstyle\,\rangle\!\rangle\ $}}
\def\qed{\nobreak\hfill $\square$ \goodbreak}

\def\limind{\mathop{\oalign{lim\cr \hidewidth$\longrightarrow$\hidewidth}}}
\def\limproj{\mathop{\oalign{lim\cr \hidewidth$\longleftarrow$\hidewidth}}}

\centerline{\bf TOPOLOGICAL GROUPS OF KAC-MOODY TYPE,}

\pec

\centerline{\bf FUCHSIAN TWINNINGS AND THEIR LATTICES}

\gec

\centerline{Bertrand R\'EMY$^*$ \quad\quad Mark RONAN$^{**}$}

\gec

\centerline{\sevenrm $^*$ UMR 5582, Institut Fourier, BP 74, 38402 
Saint-Martin-d'H\`eres Cedex, France}

\centerline{\sevenrm $^{**}$ Department of Mathematics, Statistics, 
and Computer Science,
University of Illinois at Chicago} 
\centerline{\sevenrm 322 Science and Engineering Offices (SEO) m/c 249
851 S. Morgan Street, Chicago, IL 60607-7045, USA}

\vskip 10mm

\centerline{Pr\'epublication de l'Institut Fourier n$^0$ {\bf 563} (2002)}

\pec

\centerline{\tt http://www-fourier.ujf-grenoble.fr/prepublications.html}

\vskip 10mm

A{\sevenrm BSTRACT}.--- \it
This paper deals with a class of totally disconnected groups acting 
on buildings, among
which are certain Kac-Moody groups.
The apartments of our buildings are hyperbolic planes tiled by 
right-angled polygons.
We discuss linearity properties for the groups, as well as an analogy 
with semisimple groups
over local fields of positive characteristic.
Looking for counter-examples to this analogy leads to the 
construction of Moufang twinnings
\og with several ground fields\fg.
\rm

\vskip 10mm

\footnote{}{\sevenrm {\sevenbf Keywords:} Twin buildings, twin root 
data, trees, hyperbolic buildings, lattices, non-linear
groups

}

\centerline{\bf Introduction}

\pec

J. Tits' definition of Kac-Moody groups over fields [T2] can be seen 
as the achievement of a long
process, which involved many people and many points of view.
Among previous constructions, are important works by V. Kac and D. 
Peterson [KP] and by O.
Mathieu [Mat].
Whereas O. Mathieu made use of delicate properties from algebraic 
geometry to derive
representation-theoretic results (via Schubert varieties), the 
article [KP] proposed a
combinatorial structure (refining that of a $BN$-pair) to study 
Kac-Moody groups over fields of
characteristic 0 (seen as automorphism groups of Lie algebras).
V. Kac and D. Peterson obtained some defining relations which reappeared in the
generalized Steinberg presentation used by J. Tits [T2].
Later J. Tits gave a list of group-theoretic axioms well-suited to the
Kac-Moody situation [T5], which he called {\it twin root data}; its 
geometric counterpart is the
theory of {\it twin buildings}.

\pec

The work in [T2] describes Kac-Moody groups $\Lambda$ generated by 
root groups, like
semisimple algebraic groups, except that the Weyl group $W$ is infinite.
The positive (resp. negative) root groups, along with a torus, 
generate a subgroup $B_+$ (resp.
$B_-$), but $B_+$ and $B_-$ are not conjugate because $W$ fails to 
contain an element of
longest length performing the conjugation.
As a result, there are two buildings $\Delta_+$ and $\Delta_-$, whose 
chambers can be regarded
as conjugates of $B_+$ and $B_-$ respectively.
Two chambers in the same building cannot be opposite, but chambers in 
one can be opposite
chambers in the other.
This was the starting point for the theory of twin buildings, due 
initially to J. Tits and the second
author.

\pec

As a special case, consider the Kac-Moody group $\Lambda= {\rm 
SL}_n(k[t,t^{-1}])$ over the
field $k$.
The Weyl group is of affine type $\widetilde A_{n-1}$, and the subgroups
${\rm SL}_n(k[t])$ and ${\rm SL}_n(k[t^{-1}])$ stabilize vertices in $\Delta_+$
and $\Delta_-$ respectively.
When $k$ is a finite field ${\bf F}_q$, the buildings $\Delta_+$ and 
$\Delta_-$ are locally finite, and
the facet stabilizers of $\Delta_-$ act discretely on $\Delta_+$.
On the other hand if we forget about $\Delta_-$, then in its action 
on $\Delta_+$, $\Lambda$ has a
natural topological closure $G$, namely ${\rm SL}_n \bigl( {\bf F}_q 
(\!( t )\!) \bigr)$.
Thus the twin root datum that generates $\Lambda$ gives rise to 
lattices such as
${\rm SL}_n({\bf F}_q[t^{-1}])$ inside the topological group
$G={\rm SL}_n \bigl( {\bf F}_q (\!( t )\!) \bigr)$.
This production of lattices applies to all twin root data for which 
the buildings are locally finite.
The data may or may not come from a Kac-Moody group.
If it does, we call its closure $G$ a {\it topological Kac-Moody 
group}, and if not then
a {\it topological group of Kac-Moody type}.

\pec

Q{\sevenrm UESTION}.--- \it
To what extent is a topological group of Kac-Moody type analogous to 
a semisimple group over a local field?
\rm\pec

If the group $G$ is a topological Kac-Moody group, the following theorem,
stated more precisely in 1.C.2, gives one analogy:

\pec

T{\sevenrm HEOREM}.--- \it
A $($locally compact, totally disconnected$)$ topological Kac-Moody 
group admits a refined Tits system.
This Tits system gives rise to a building in which any 
$($spherical$)$ facet-fixator is, up to finite index,
a pro-$p$ group.
\rm\pec

A {\it refined Tits system} is the structure proposed by V. Kac and D.
Peterson to study Kac-Moody groups $\Lambda$ [KP].
It is defined in 1.A.5, and unlike a twin root datum, defined in 
1.A.1, is asymmetric with respect to
positive and negative root groups, because it is concerned with only 
one building rather than two.
This building is identified with its nonpositively curved realization 
[D], where only
spherical residues appear as facets.
The fixators of those facets are analogues of the parahoric subgroups 
of Bruhat-Tits theory, and their
finite index pro-$p$ subgroups are analogues of congruence subgroups.

\pec

The above analogy between semisimple groups over local fields and 
topological Kac-Moody groups can
be extended to the more general setting of topological groups of 
Kac-Moody type.
We do this by constructing twin root data that do not arise from 
Kac-Moody groups.
This gives exotic twin buildings.
In view of the classification work [M\"u, MR] on 2-spherical twin 
buildings (meaning that their rank 2 residues have finite
(dihedral) Weyl groups), the twin buildings we construct are not 2-spherical.
One consequence of our work is the existence of Moufang twin 
buildings (i.e. twin buildings admitting root groups), of any
desired finite rank, that do not arise from Kac-Moody groups.
The apartments in these buildings can be realised as tilings of the 
hyperbolic plane by right-angled $r$-gons
($r \geq 5$).
The rank of such a building is $r$, and the Weyl group is generated 
by reflections across the faces of a fundamental domain $R$
in the hyperbolic plane.
Buildings with such apartments are called {\it Fuchsian } by M. 
Bourdon [Bou3], and he proves existence and uniqueness
results for them.
Uniqueness [loc.cit.], along with the construction given in the 
present paper, proves the following theorem,
described in more detail in 4.E.2:

\pec

T{\sevenrm HEOREM}.--- \it
A right-angled Fuchsian building belongs to a Moufang twinning 
whenever its thicknesses at panels
are cardinalities of projective lines.
\rm\pec

The group combinatorics obtained by standard arguments [A] from these 
buildings are not Kac-Moody
as soon as we choose two thicknesses \og of different characteristics\fg.
As already mentioned, these buildings are non-2-spherical, since the 
reflections across any two
non-intersecting edges of the polygon $R$ generate an infinite 
dihedral group. Hence, they are not covered
by the Local to Global theorem [MR] that applies in the 2-spherical 
case, and is leading
to a classification [M\"u].
On the other hand, techniques from CAT$(-1)$-geometry lead to deep
rigidity and amenability theorems for (lattices of) Fuchsian 
buildings [BM, BP].
 From the point of view of lattices of CAT$(-1)$ geometries, we can 
formulate the following consequence
of our construction (5.B).

\pec

T{\sevenrm HEOREM}.--- \it
There exist groups with twin root data such that the associated buildings are
Fuchsian, and in which any Borel subgroup of a given sign is a 
non-uniform lattice of the building of opposite sign.
Moreover any group homomorphism from such a lattice to a product of 
linear algebraic groups $($possibly over
different fields$)$ has infinite kernel.
\rm\pec

In particular, the above groups with twin root data are not linear 
over any field.
This shows that the remaining question concerning linearity of 
Kac-Moody groups -- to find a
Kac-Moody group not linear over any field [R\'e3] -- is solved in the wider
(easier) context of groups with twin root data.
Discussing linearity, even in the smaller class of Kac-Moody groups, 
leads to surprising
situations, as shown by our last result (5.C).

\pec

P{\sevenrm ROPOSITION}.--- \it
There exist topological Kac-Moody groups over ${\bf F}_q$ which admit both
non-uniform lattices which cannot be linear over any field of 
characteristic prime to $q$,
and uniform lattices which have convex-cocompact embeddings into real 
hyperbolic spaces.
\rm\pec

The existence of the latter lattices and of their linear embeddings 
is due to M. Bourdon [Bou1], who
proves that the limit sets of the embeddings often have Hausdorff 
dimension $>2$.

\pec

This paper is organized as follows.
In the first section, we recall some combinatorial facts on twin root 
data (1.A).
Then we define the topological groups of Kac-Moody type and 
topological Kac-Moody groups (1.B).
Subsect. 1.C provides arguments to see the latter class of groups as 
a generalization of
algebraic groups over local fields of positive characteristic.
Sect. 2 sketches a construction of twin root data, and recalls some 
relations in SL$_2$.
Sect. 3 applies this construction to the case of trees.
It is a special case of the construction of Moufang twin trees due to 
J. Tits [T4], but given in a very
down-to-earth way.
The details we supply in Sect. 3 allow us to avoid computation in 
Sect. 4 where we concentrate on
the geometry of Fuchsian buildings.
Finally, in Sect. 5 we discuss linearity properties for lattices of topological
groups of Kac-Moody type: this is relevant to the study of Kac-Moody 
lattices as generalizations of
arithmetic groups over function fields [R\'e3].

\pec

Let us state a convention we will use for group actions.
If a group $G$ acts on a set $X$, the (pointwise) stabilizer of a 
subset $Y \subset X$ will
be called its {\it fixator}, and will be denoted by ${\rm Fix}_G(Y)$.
The classical (global) stabilizer will be denoted by ${\rm Stab}_G(Y)$.
When $Y$ is a facet of a building $X$, and $G$ is a type preserving group of
automorphisms, the two notions of course coincide.
Finally the notation $G \! \mid_X$ means the group obtained from $G$ 
by factoring out the
kernel of the action on $X$.

\pec

The first author thanks University College of London where this work 
was initiated, the Hebrew
University for its warm hospitality during the academic year 
2000/2001, and more personally,
Shahar Mozes for suggesting a strong version of the non-linearity 
property 5.B, and Marc Bourdon (resp.
Guy Rousseau) for helpful discussions about 5.C (resp. 1.C).
The first author was partially supported by grants from the British 
EPSRC and the Hebrew University; the second
author was partially supported by the National Security Agency.

\vfill\eject

\centerline{\bf 1. Closures of groups with twin root data}

\pec

The purpose of this section is to introduce a family of totally disconnected
topological groups, called {\it topological groups of Kac-Moody type} 
-- see 1.B.
This requires some basic properties of group combinatorics introduced 
by J. Tits
and by V. Kac and D. Peterson, which we describe in Subsect. 1.A.
In the case of topological Kac-Moody groups, some precise structure results are
available, which are proved in 1.C.
We exhibit analogies with an algebraic group over a local field
of characteristic $p$, as mentioned above.

\pec

{\bf 1.A } {\it Combinatorial framework}. --- In this subsection, we 
present all the combinatorial material we will need in the sequel.
For references, see for instance [A, KP, R\'e1, Ro1-2, T2-5].

\pec

{\bf 1.A.1 } Let $(W,S)$ be the Coxeter system

\pec

\centerline{$W = \langle s \!\in \! S \mid (st)^{M_{st}}=1$ whenever 
$M_{st}<\infty
\rangle$,}

\pec

where $M$ is a Coxeter matrix.
The existence of an abstract simplicial complex $\Sigma=\Sigma(W,S)$ 
acted upon by $W$ is the
starting point for the definition of buildings of type $(W,S)$ in 
terms of apartment systems.
This complex is the {\it Coxeter complex } associated to $W$ [Ro1 \S2].
Set-theoretically, $\Sigma$ is the union of the translates
$wW_J:= w\langle J \rangle$ for $J \subset S$ and $w \! \in \! W$.
It is partially ordered by reverse inclusion.
The elements $w \langle \varnothing \rangle$, which are simply the 
elements of $W$, have
maximum dimension and are called {\it chambers}.
The {\it root system } of $(W,S)$ is defined by means of the length function
$\ell: W \to {\bf N}$ with respect to $S$ [T2 \S5].
The set $W$ admits a $W$-action via left translations.
{\it Roots } are distinguished halves of $W$, regarded as a $W$-set.

\pec

D{\sevenrm EFINITION}.--- \it
{\rm (i)} The {\rm simple root } of index $s \! \in \! S$ is the half
$a_s:=\{ w \! \in \! W \mid \ell(sw) > \ell(w) \}$.

{\rm (ii) } A {\rm root} of $W$ is a translate $wa_s$, $w \! \in \! W,
s \! \in \! S$.
The set of all roots will be denoted by $\Phi$.
It admits an obvious $W$-action by left translations.

{\rm (iii) } A root is {\rm positive } if it contains
$1$; otherwise, it is {\rm negative}.
We denote by $\Phi_+$ $($resp. $\Phi_-)$ the set of positive $($resp. 
negative$)$ roots.

{\rm (iv) } The complement of a root $a$, denoted $-a$, is also a 
root, called the
{\rm opposite} of $a$.

{\rm (v) }The {\rm boundary }$($or {\rm wall}$)$ $\partial a$ of a 
root $a$ is the union of the closed panels
having exactly one chamber in $a$.
This is the same as the intersection of the closures of $a$ and $-a$.
\rm\pec

The following definitions are needed for the group combinatorics.

\pec

D{\sevenrm EFINITION}.--- \it
{\rm (i) } A pair of roots $\{ a;b \}$ is called {\rm prenilpotent } if both
intersections $a \cap b$ and $(-a) \cap (-b)$ are
nonempty.

{\rm (ii) } Given a prenilpotent pair of roots $\{a; b\}$,
the {\rm interval }$[a; b]$ is by definition the set of roots
$c$ with $c \supset a \cap b$ and
$(-c) \supset (-a) \cap (-b)$.
\rm\pec

We can now turn to groups.

\pec

D{\sevenrm EFINITION}.--- \it
Let $\Lambda$ be an abstract group containing a subgroup $H$.
Suppose $\Lambda$ is endowed with a family $\{ U_a \}_{a \in \Phi}$ 
of subgroups
indexed by the set of roots $\Phi$, and define the subgroups
$U_+:= \langle U_a \mid a \! \in \! \Phi_+ \rangle$ and
$U_-:= \langle U_a \mid a \! \in \! \Phi_- \rangle$.
We say that the triple $\bigl( \Lambda, \{ U_a \}_{a \in \Phi}, H \bigr)$
is a {\rm twin root datum } for $\Lambda$ $($or satisfies the {\rm 
(TRD) axioms}$)$ if the
following conditions are fulfilled.

\pec

{\rm (TRD0) } Each $U_a$ is nontrivial and normalized by $H$.

{\rm (TRD1) } For each prenilpotent pair of roots
$\{ a; b \}$, the commutator subgroup $[U_a, U_b]$ is
contained in the subgroup $\langle U_c: c \! \in [a; b]-\{a;b\} \rangle$.

{\rm (TRD2) } For each $s \! \in \! S$ and $u \! \in \! U_{a_s} 
\setminus \{ 1 \}$,
there exist uniquely defined $u', u'' \! \in \! U_{-a_s} \setminus \{ 
1 \}$ such that
$m(u):=u'uu''$ conjugates $U_b$ onto $U_{s.b}$ for every root $b$.
Moreover $m(u)H=m(v)H$ for all $u, v \! \in \! U_{a_s} \setminus \{ 1 \}$.

{\rm (TRD3) } For each $s \! \in \! S$, we have
$U_{a_s} \not\subset U_-$ and $U_{-a_s} \not\subset U_+$.

{\rm (TRD4) } The group $\Lambda$ is generated by $H$ and the $U_a$'s.
\rm\pec

E{\sevenrm XAMPLE}. --- The special linear groups ${\rm SL}_n({\bf F}_q[t,t^{-1}])$
over Laurent polynomials satisfy these axioms for $W$ affine of type 
$\widetilde A_{n-1}$.
More generally, the groups ${\Bbb G}({\bf F}_q[t,t^{-1}])$ where 
${\Bbb G}$ is a
simply connected Chevalley group, satisfy these axioms when $W$ is a 
suitable affine Coxeter group [A
\S1]. In these cases, the root groups are naturally defined as 
subgroups of the spherical root groups by
using polynomials of the form
$at^r$ for fixed $r$ as the off-diagonal entries.

\pec

{\bf 1.A.2 }
Here are the main properties of such a group $\Lambda$.
Define the {\it standard Borel subgroup of sign $\epsilon=\pm$ } to 
be $HU_\epsilon$.
The subgroup $N<\Lambda$ is by definition generated by $H$ and the
$m(u)$'s of axiom (TRD2). Then $(\Lambda,HU_+,N,S)$ and 
$(\Lambda,HU_-,N,S)$ are $BN$-pairs
sharing the same Weyl group $W=N/H$ [T3-4], and $HU_+$ and $HU_-$ are
conjugate if and only if $W$ is finite.

\pec

A formal consequence of the existence of a $BN$-pair is a Bruhat 
decomposition for the
group, and using the root groups above this can be done as follows.
For each $w \! \in \! W$, define the subgroups
$U_w:= U_+ \cap wU_-w^{-1}$ and $U_{-w}:= U_- \cap wU_+w^{-1}$.
In their action on the buildings defined by the $BN$-pairs above, the 
groups $U_w$ and $U_{-w}$ are freely transitive on the
sets of chambers at distance $w$ from the (base) chambers fixed by 
$HU_+$ and $HU_-$ respectively.

\pec

For each $w \! \in \! W$, define the (finite) sets of roots

\pec

\centerline{
$\Phi_w:= \Phi_+ \cap w^{-1}\Phi_-$ and $\Phi_{-w}:= \Phi_- \cap 
w^{-1}\Phi_+$.}

\pec

In terms of buildings, if $c_+$ represents the chamber fixed by 
$HU_+$ in the apartment stabilised by $N$, then $\Phi_w$ is
the set of roots in that apartment containing $c_+$ but not 
containing $w^{-1}c_+$.
Similarly for $\Phi_{-w}$ with respect to the chamber fixed by $HU_-$.
The group $U_w$ (resp. $U_{-w}$) is in bijection with the set-theoretic
product of the root groups indexed by $\Phi_{w^{-1}}$ (resp. 
$\Phi_{-w^{-1}}$) for a suitable
(cyclic) ordering [T4, proposition 3].
Note that $U_w<U_{w'}$ as soon as $w \leq w'$ for the Bruhat 
ordering, because $\Phi_{w^{-1}}$ is a
subset of $\Phi_{w'^{-1}}$ in this case.

\pec

The {\it refined Bruhat decompositions } are [KP, Proposition 3.2]:

\pec

\centerline{$\Lambda = \bigsqcup_{w \in W} U_wwHU_+$ \quad and \quad
$\Lambda = \bigsqcup_{w \in W} U_{-w}wHU_-$.}

\pec

In these decompositions, the element in the left factor $U_{\pm w}$ 
is uniquely determined.

\pec

A third kind of decomposition involves both signs and is used to 
define the so-called {\it twin structures} (twin buildings, twin
$BN$-pairs... ).
The {\it refined Birkhoff decompositions } are [KP, Proposition 3.3]:

\pec

\centerline{$\Lambda = \bigsqcup_{w \in W} (U_+ \cap wU_+w^{-1})wHU_-
= \bigsqcup_{w \in W} (U_- \cap wU_-w^{-1})wHU_+$.}

\pec

{\bf 1.A.3 } The theory of {\it $($Moufang$)$ twin buildings } is the 
geometric side of the
group combinatorics above [T3-T5].
The definition of a twin building is quite similar to that of a 
building in terms of $W$-distance
[T5, A \S 2, Ro2 \S 2].

\pec

D{\sevenrm EFINITION}.--- \it
A {\rm twin building of type $(W,S)$ } consists of two buildings 
$(\Delta_+,d_+)$ and $(\Delta_-,d_-)$ of type $(W,S)$
endowed with a {\rm ($W$-)codistance}, namely a map

\pec

\centerline{
$d^*: (\Delta_+\times\Delta_-) \cup (\Delta_-\times\Delta_+) \rightarrow W$}

\pec

satisfying the following conditions for each sign
$\epsilon$ and all chambers $x_\epsilon$ in $\Delta_\epsilon$ and
$y_{-\epsilon}$, $y_{-\epsilon}'$ in $\Delta_{-\epsilon}$.

\pec

{\rm (TW1) } 
$d^*(y_{-\epsilon},x_\epsilon)=d^*(x_\epsilon,y_{-\epsilon})^{-1}$.

{\rm (TW2) } If $d^*(x_\epsilon,y_{-\epsilon})=w$ and
$d_{-\epsilon}(y_{-\epsilon},y_{-\epsilon}')=s \! \in \! S$ with 
$\ell(ws)<\ell(w)$,
then $d^*(x_\epsilon,y_{-\epsilon}')=ws$.

{\rm (TW3) } If $d^*(x_\epsilon,y_{-\epsilon})=w$
then for each $s \! \in \! S$,
there exists $z_{-\epsilon} \! \in \! \Delta_{-\epsilon}$ with
$d_{-\epsilon}(y_{-\epsilon},z_{-\epsilon})=s$ and
$d^*(x_\epsilon,z_{-\epsilon})=ws$.
\rm\pec

Two chambers are {\it opposite } if they are at codistance $1$; two 
facets are {\it
opposite } if they belong to opposite chambers and have the same type.
Given an apartment ${\Bbb A}_\epsilon$ of sign $\epsilon$, an {\it opposite }
of it is an apartment
${\Bbb A}_{-\epsilon}$ such that each chamber of ${\Bbb A}_\epsilon$ 
is opposite exactly one chamber in ${\Bbb
A}_{-\epsilon}$.
An apartment has at most one opposite, and a pair of opposite 
apartments is called a {\it twin
apartment} [A \S 2, Ro2, T3-5].
Note that every pair of opposite chambers lies in exactly one twin apartment.

\pec

Set-theoretically, the set of chambers of the building of sign 
$\epsilon$ attached to a twin root datum
$\bigl( \Lambda, \{ U_a \}_{a \in \Phi}, H \bigr)$ is
$\Lambda/HU_\epsilon=\{gHU_\epsilon\}_{g \in \Lambda}$.
The $\Lambda$-action is the natural one by left translations, and it 
is transitive on the pairs of opposite
chambers -- see [A \S2] for details.

\pec

E{\sevenrm XAMPLE}. --- When $\Lambda$ is ${\Bbb G}({\bf 
F}_q[t,t^{-1}])$ where ${\Bbb G}$ is a
Chevalley group, the buildings $\Delta_+$ and $\Delta_-$ are the 
Bruhat-Tits buildings of
${\Bbb G} \bigl({\bf F}_q (\!( t )\!) \bigr)$ and ${\Bbb G}\bigl({\bf 
F}_q (\!( t^{-1} )\!) \bigr)$
respectively, on which $\Lambda$ acts diagonally.

\pec

{\bf 1.A.4 } We state now the generalized Levi decompositions of 
parabolic subgroups (which are by
definition the fixators of facets).
Let $F$ be a spherical facet.
Then for every choice of a twin apartment ${\Bbb A}$ containing $F$, 
the parabolic subgroup ${\rm Fix}_\Lambda(F)$ admits
a semi-direct product decomposition [R\'e1, 6.2]:

\pec

\centerline{${\rm Fix}_\Lambda(F) = M(F,{\Bbb A}) \ltimes U(F)$.}

\pec

The group $M(F,{\Bbb A})$ is the fixator of $F \cup -F$, where $-F$ 
is the opposite of $F$ in
${\Bbb A}$; it is generated by $H = {\rm Fix}_\Lambda({\Bbb A})$ and
the root groups $U_a$ with $\partial a \supset F$.
Moreover it satisfies the {\rm (TRD) } axioms for the finite root 
system attached to $F$.
The group $U(F)$ only depends on $F$, and is the normal closure in 
${\rm Fix}_\Lambda(F)$ of
the root groups $U_a$ for which $a\supset F$.

\pec

{\bf 1.A.5 } We end this subsection by introducing some group 
combinatorics defined by V. Kac
and D. Peterson, which we will need for the definition of topological 
groups of Kac-Moody
type.

\pec

D{\sevenrm EFINITION}.--- \it
A {\rm refined Tits system } for a group $G$ is a sixtuple
$(G,N,U_+,U_-,H,S)$, where $N$, $U_+$, $U_-$ are subgroups satisfying 
the following
axioms:

\pec

{\rm (RT1)} We have $G = \langle N, U_+ \rangle $, $H \triangleleft N$,
$H<N_G(U_+) \cap N_G(U_-)$, $W:=N/H$ and $(W,S)$ is a Coxeter system.

{\rm (RT2)} For each $s \! \in \! S$, we set $U_s:=U_+ \cap s^{-1} 
U_- s$; and for any
$w \! \in \! W$ and $s \! \in \! S$ we require:

\hskip 5mm {\rm (RT2a)} \quad $s^{-1} U_s s \neq \{ 1 \}$ and
$s^{-1}( U_s \setminus \{ 1 \} ) s \subset U_s sH U_s$,

\hskip 5mm {\rm (RT2b)} \quad either $w^{-1} U_s w \subset U_+$ or 
$w^{-1} U_s w \subset U_-$,

\hskip 5mm {\rm (RT2c)} \quad $U_+ = U_s(U_+ \cap s^{-1} U_+ s)$.

{\rm (RT3)} If $u_-\! \in \! U_-$, $u_+\! \in \! U_+$ and $n\! \in \! 
N$ are such that
$u_-nu_+=1$, then $u_-=u_+=n=1$.
\rm\pec

Many properties can be derived from these axioms: refined Bruhat 
decomposition, definition of
subgroups as limits of natural inductive systems, etc.
The reference is the original paper [KP].
We will see in 1.C why this framework concerns Lie groups over 
non-Archimedean local fields of
{\it positive } characteristic, but not the characteristic 0 case.
The relation between twin root data and refined Tits sytems is the 
following -- see [R\'e1, 1.5.4] for a proof.

\pec

T{\sevenrm HEOREM}.--- \it
If $\bigl( \Lambda, \{U_a \}_{a \in \Phi}, H \bigr)$ is a twin root 
datum indexed by the Coxeter system $(W,S)$, then
$(\Lambda,N,U_+,U_-,H,S)$ and $(\Lambda,N,U_-,U_+,H,S)$ are refined 
Tits system with Weyl group $W$.
Moreover we have $U_+ \cap s^{-1} U_- s=U_{a_s}$ and $U_- \cap s^{-1} 
U_+ s=U_{-a_s}$
for any $s \! \in \! S$.
\qed\rm\pec

{\bf 1.B } {\it Topological groups of Kac-Moody type}. ---
In this subsection, we use the buildings of groups with twin root 
data to define related topological groups.
This requires to see buildings as metric spaces.

\pec

{\bf 1.B.1 } We assume we are given a group $\Lambda$ endowed with a 
twin root datum
$\bigl( \Lambda, \{ U_a \}_{a \in \Phi}, H \bigr)$.
We only deal with the positive building of $\Lambda$, which we simply 
denote by $\Delta$.
The notation $\Delta$ actually denotes the metric realization of the 
building in the sense of
Moussong-Davis [D, Mou]. Let us list some assumptions and discuss 
their consequences.

\pec

A{\sevenrm SSUMPTION} 1.--- \it
All root groups $U_a$ are finite.
\rm\pec

We can thus set

\pec

\centerline{$\displaystyle q_{\hbox{\sevenrm min}}:=
\min_{a \in \Phi} \mid \! U_a \! \mid = \min_{s \in S} \mid \! U_{a_s} \! \mid$
\quad and \quad
$\displaystyle q_{\hbox{\sevenrm max}}:=
\max_{a \in \Phi} \mid \! U_a \! \mid = \max_{s \in S} \mid \! 
U_{a_s} \! \mid$.}

\pec

By the Moufang property [Ro1 \S6], the number of chambers on each 
panel is then bounded by
$1 + q_{\hbox{\sevenrm max}}$.
More generally, the fact that only spherical facets appear in the 
Moussong-Davis realization of a building implies that
$\Delta$ is a locally finite cell complex.

\pec

A{\sevenrm SSUMPTION} 2.--- \it
The Weyl group $W$ is infinite; hence so is the root system $\Phi$.
\rm\pec

This implies by [D] that the building $\Delta$ is a non-positively 
curved metric space, which
means that it satisfies the so-called CAT(0)-property [BH, II.1.1]: 
geodesic triangles are at least as thin
as in the Euclidean plane.

\pec

Our main motivation in defining topological groups of Kac-Moody type 
is geometric group theory, in which
groups are studied via their actions on suitable spaces.
The following lemma shows that moding out by the kernel of the 
$\Lambda$-action on $\Delta$ is
harmless for the group combinatorics.

\pec

L{\sevenrm EMMA}.--- \it
Let $\bigl( \Lambda, \{U_a \}_{a \in \Phi}, H \bigr)$ be a twin root 
datum with associated Coxeter system $(W,S)$ and
associated buildings $\Delta_\pm$.
Then the kernel $K$ of the $\Lambda$-action on $\Delta_+$ or 
$\Delta_-$ consists of the elements in $H$
centralizing all root groups.
The groups $U_\pm$ embed in $\Lambda/K$, hence so do the root groups, and
$\bigl( \Lambda/K, \{U_a \}_{a \in \Phi}, H/K \bigr)$ is again a twin 
root datum with the same
associated Coxeter system and twinning.
\rm\pec

{\it Proof}. We use theorem 1.A.5, which allows us to use properties 
of refined Tits systems for $\Lambda$.
We argue with $\Delta_+$, the negative case being completely similar.
The kernel $K$ fixes the chamber fixed by $HU_+$ and the apartment 
stabilized by $N$; hence it lies in $HU_+ \cap N$.
According to [KP, lemma 3.1], the latter intersection is $H$.
That $K$ must centralize each root group comes from the Moufang 
property, since each root group is in
bijective correspondence with the set of chambers sharing a panel 
with a given chamber [Ro1, 6.4].
By (RT3) we have $U_\pm \cap H =\{ 1 \}$, hence $U_+$ and $U_-$ embed 
in $\Lambda/K$, which
implies axioms (TRD0) and (TRD1) (a prenilpotent pair of roots can 
always be transformed into a pair of
positive roots by a suitable element of $W$).
Axiom (TRD4) is clear and (TRD3) for $\Lambda/K$ is a consequence of 
(RT3) for $\Lambda$.
Axiom (TRD2) for $\Lambda/K$ is a consequence of axiom (TRD2) for 
$\Lambda$, the uniqueness part coming from
the fact that, when all the other axioms are satisfied, the elements 
$u'$ and $u''$ in (TRD2) are uniquely determined by $u$ [T5, p.257].
The set-theoretical definition of the positive (resp. negative) 
building as $\Lambda/HU_+$ (resp. $\Lambda/HU_-$) and the
definition of the $W$-distances and of the codistance by means of the 
Bruhat and Birkhoff decompostions show that the
associated twinnings are the same.
\qed\pec

R{\sevenrm EMARK}. ---
When $\Lambda$ is a Kac-Moody group, the kernel $K$ is the center 
$Z(\Lambda)$ [R\'e1, 9.6.2], which is finite when the
ground field is.

\pec

E{\sevenrm XAMPLE}. --- For ${\rm SL}_n({\bf K}[t,t^{-1}])$, the 
buildings of $\Lambda$ are Euclidean
$\widetilde A_{n-1}$-buildings, and the kernel $K$ is the group 
$\mu_n({\bf F}_q)$ of $n$-th roots of unity in ${\bf F}_q$.

\pec

The lemma leads us to make the following convention.

\pec

C{\sevenrm ONVENTION}.--- \it
Until the end of {\rm 1.B}, $\Lambda$ denotes the image of a group 
with twin root datum as above under the group
homomorphism corresponding to the action on the building $\Delta$.
Consequently, $\Lambda<{\rm Aut}(\Delta)$ admits a twin root datum 
with finite root groups and infinite
Weyl group.
\rm\pec

{\bf 1.B.2 } We can now define topological groups of Kac-Moody type.
Given an automorphism $g \! \in \! {\rm Aut}(\Delta)$ and a finite 
union $C$ of facets, the
subset of ${\rm Aut}(\Delta)$

\pec

\centerline{$O_C (g):= \{ h \! \in \! {\rm Aut}(\Delta): h \! \mid_C 
= g \! \mid_C \}$}

\pec

is by definition an open neighborhood of $g$.
The group $\Lambda$ is not discrete on the single building $\Delta$ 
because each facet
fixator is transitive on opposite facets in $\Delta_-$, of which 
there are infinitely many.

\pec

D{\sevenrm EFINITION}.--- \it
Let $\Lambda$ be a group as in convention {\rm 1.B.1}, with finite
root groups and infinite Weyl group.

\pec

{\rm (i) } We call the closure

\pec

\centerline{
$G:= \overline\Lambda^{{\hbox{\sevenrm Aut}}(\Delta)}$}

\pec
a {\rm topological group of Kac-Moody type}.
If $\Lambda$ is a Kac-Moody group over a finite field, we call $G$ a 
{\rm topological Kac-Moody group}.

{\rm (ii) } We call the fixator ${\rm Fix}_G(F)$ of a facet $F$ the 
{\rm parahoric subgroup} of $G$ associated to $F$,
and we denote it by $G_F$.
We call a chamber fixator an {\rm Iwahori subgroup}.

{\rm (iii) } We denote by $U_F$ the closure 
$\overline{U(F)}^{{\hbox{\sevenrm Aut}}(\Delta)}$ of the
\og unipotent radical\fg of ${\rm Fix}_\Lambda(F)$.
\rm\pec

R{\sevenrm EMARK}. ---
This definition implies that an automorphism of $\Delta$ lies in $G$ 
if and only if it coincides
with an element of $\Lambda$ on any finite set of facets.

\pec

E{\sevenrm XAMPLE}. --- The image of ${\rm SL}_n({\bf 
F}_q[t,t^{-1}])$ under the action on its positive Euclidean building 
is
the group ${\rm SL}_n({\bf F}_q[t,t^{-1}])/\mu_n({\bf F}_q)$ -- see 1.B.1.
Then the \og completion\fg $G$ is the ultrametric Lie group
${\rm SL}_n \bigl( {\bf F}_q(\!(t)\!) \bigr) /\mu_n({\bf F}_q)$.
This is close to being ${\rm Aut}(\Delta)$ when $n>2$, but when 
$n=2$, $\Delta$ is a tree and $G$ is far from being all of
${\rm Aut}(\Delta)$.
(Locally, the difference appears in the action on vertex stars: in $G$ one has
${\rm PSL}_2({\bf F}_q)$, but in ${\rm Aut}(\Delta)$ one has the 
symmetric group ${\cal S}_{q+1}$.)

\pec

The extension of the classical terminology of parahoric subgroups, 
supported by the above example, suggests an analogy with
algebraic groups over local fields.
One of the main questions in this work is to understand to what 
extent the analogy is relevant to this more general setting.

\pec

{\bf 1.B.3 } We now consider discrete subgroups of the topological group $G$.
This will make the analogy deeper, bringing in a certain class of lattices
generalizing some arithmetic groups over function fields.
Fix a chamber $c_-$ in $\Delta_-$, and let $\Gamma$ denote its stabilizer.
Let ${\Bbb A}$ be a twin apartment containing $c_-$.
The combinatorial properties derived from the (TRD) axioms prove the 
following results:

\pec

-- The positive apartment ${\Bbb A}_+$ is a fundamental domain for 
the $\Gamma$-action on $\Delta$.
More generally, if $F_-$ is a facet in ${\Bbb A}_-$, we denote by 
$F_+$ its unique opposite in
${\Bbb A}_+$.
Given a chamber $c_+$ containing $F_+$ in its closure, each wall of 
${\Bbb A}_+$
containing $F_+$ bounds a root containing $c_+$.
The intersection of these half-spaces is a fundamental domain for the action of
${\rm Fix}_\Lambda(F_-)$ on $\Delta$. See [A \S3] for details.

-- The group $\Lambda$ acts diagonally as a discrete group on the 
product of buildings
$\Delta_+ \times \Delta_-$.
Consequently, the group $\Gamma$, hence the fixator (in $\Lambda$) of 
any negative facet is a discrete
subgroup of $G$ [CG, R\'e2].

-- If $q_{\hbox{\sevenrm min}}$ is large enough, then the (spherical) 
parabolic subgroups of a
given sign $\epsilon$ are lattices of the locally finite 
CAT(0)-building $\Delta_{-\epsilon}$
[CG, R\'e2].

\pec

In the same spirit, we can also prove a result about commensurators.
Recall that two subgroups $\Gamma$ and $\Gamma'$ of a given group $G$ 
are {\it commensurable } if
they share a finite index subgroup.
The {\it commensurator} of $\Gamma$ in $G$ is the group

\pec

\centerline{${\rm Comm}_G(\Gamma):= \{ g \! \in \! G: \Gamma$ and
$g \Gamma g^{-1}$ are commensurable$\}$.}

\pec

L{\sevenrm EMMA}.--- \it
Let $F$ be a negative $($spherical$)$ facet.

\pec
{\rm (i) } For any other negative $($spherical$)$ facet $F'$, the fixators
${\rm Fix}_\Lambda(F)$ and ${\rm Fix}_\Lambda(F')$ are commensurable.

{\rm (ii) } The group $\Lambda$ lies in the commensurator
${\rm Comm}_{{\hbox{\sevenrm Aut}}(\Delta)} \bigl( {\rm Fix}_\Lambda(F) \bigr)$
of the facet fixator ${\rm Fix}_\Lambda(F)$.
\rm\pec

{\it Proof}. We keep our negative Borel subgroup $\Gamma = {\rm 
Fix}_\Lambda(c_-)$.
The arguments for both points are based on finiteness of root groups 
and refined Bruhat decomposition.
Recall that in 1.A.2 we attached a finite set of roots $\Phi_w$ to 
any element of the Weyl group
$w \! \in \! W$.
According to [KP, proposition 3.2] and [T4, proposition 3], for any 
$z \! \in \! W$ we have:

\pec

\centerline{$\displaystyle (\star) \quad
\Gamma = \bigl( \prod_{a \in \Phi_{z^{-1}}} U_{-a} \bigr)
\cdot \bigl( \Gamma \cap z \Gamma z^{-1} \bigr)$
\quad hence \quad
$\displaystyle (\star\star) \quad
\Gamma z \Gamma = \bigl( \prod_{a \in \Phi_{z^{-1}}} U_{-a} \bigr) z \Gamma$,}

\pec

where the product of the root groups is $\{ 1 \}$ if $z=1$ (in which 
case $\Phi_{z^{-1}}$ is
the empty set).
Notice that in both equalities the first factor on the right is finite.

\pec

(i). For any spherical facet $F$, choose a chamber whose closure contains it.
By the Bruhat decomposition for parabolics and $(\star\star)$, the 
Borel subgroup
corresponding to this chamber is of finite index in ${\rm Fix}_\Lambda(F)$.
Hence we are reduced to proving that any two Borel subgroups are commensurable.
By transitivity of $\Lambda$ on pairs of chambers at given 
$W$-distance, it is enough
to consider $\Gamma$ and $w\Gamma w^{-1}$ for $w \! \in \! W$.
Setting $z=w$ in equality $(\star)$ shows that $\Gamma \cap w\Gamma 
w^{-1}$ is of
finite index in $\Gamma$.
Similarly setting $z=w^{-1}$, and conjugating the result by $w$ shows 
that $\Gamma \cap
w\Gamma w^{-1}$ is of finite index in $w\Gamma w^{-1}$.

\pec

(ii). Since $\Lambda$ is transitive on chambers and the fixator of 
any spherical facet
contains a Borel subgroup of finite index, it is enough to show that 
$\Lambda$ lies in the
commensurator of $\Gamma$. Let $g \! \in \! \Lambda$.
By $(\star\star)$, $g=u_{-w}w\gamma$ for some $w \! \in \! W$, for 
$\gamma \! \in \!
\Gamma$ and for $u_{-w}$ in the (finite) group $U_{-w}$ defined in 1.A.2.
We have then $g\Gamma g^{-1} \cap \Gamma = u_{-w}(w\Gamma w^{-1} \cap
\Gamma)u_{-w}^{-1}$ because $u_{-w} \! \in \! \Gamma$.
By (i), $w \Gamma w^{-1} \cap \Gamma$ is of finite index in $\Gamma$, and
since $u_{-w} \! \in \! \Gamma$, we have $[\Gamma: g\Gamma g^{-1} 
\cap \Gamma]<\infty$.
In other words, if $g \! \in \! \Lambda$, then $g\Gamma g^{-1} \cap 
\Gamma$ is of finite index in
$\Gamma$.
Replacing $g$ by $g^{-1}$and conjugating by $g$ shows that $g\Gamma
g^{-1} \cap \Gamma$ is of finite index in $g\Gamma g^{-1}$.
\qed\pec

E{\sevenrm XAMPLE}. --- Let us consider again the example of the group
${\rm SL}_n({\bf F}_q[t,t^{-1}])/\mu_n({\bf F}_q)$.
Its closures in the automorphism groups of the positive and negative 
builidings are respectively
${\rm SL}_n\bigl( {\bf F}_q(\!(t)\!) \bigr)/\mu_n({\bf F}_q)$ and
${\rm SL}_n\bigl( {\bf F}_q(\!(t^{-1})\!) \bigr)/\mu_n({\bf F}_q)$.
A negative vertex fixator in ${\rm SL}_n\bigl( {\bf 
F}_q(\!(t^{-1})\!) \bigr)/\mu_n({\bf F}_q)$ is isomorphic to
${\rm SL}_n \bigl( {\bf F}_q[[ t^{-1} ]] \bigr)/\mu_n({\bf F}_q)$; 
consequently, the lattices we get by taking
negative facet fixators in $\Lambda$ are $\{ 0 \}$-arithmetic groups 
commensurable with
${\rm SL}_n({\bf F}_q[t^{-1}])/\mu_n({\bf F}_q)$.
The commensurator of ${\rm SL}_n({\bf F}_q[t^{-1}])/\mu_n({\bf F}_q)$ contains
${\rm SL}_n \bigl( {\bf F}_q(t) \bigr)/\mu_n({\bf F}_q)$.

\pec

{\bf 1.C } {\it The special case of Kac-Moody groups}. ---
Let $\Lambda$ be an almost split Kac-Moody group.
According to J. Tits, a split Kac-Moody group is defined by 
generators and Steinberg
relations, once a generalized Cartan matrix $A$ and a ground field 
${\bf K}$ are fixed [T2] -- see
also 2.E.1.
An almost split Kac-Moody group is the group of fixed points of a 
suitable Galois action on a split
group -- see [R\'e1, 11.3] for a precise definition.
The basic result we use about almost split Kac-Moody groups is that 
their rational points satisfy the
(TRD) axioms [R\'e1, 12.6.3].

\pec

{\bf 1.C.1 } We make the following convention until the end of the section.

\pec

C{\sevenrm ONVENTION}.--- \it
The group $\Lambda$ is the image of the rational points of an almost 
split Kac-Moody group over
the ground field ${\bf F}_q$ $($of cardinality $q$ and characteristic 
$p)$ in the full automorphism group
of the positive building $\Delta$.
We denote its twin root datum by $\bigl( \Lambda, \{ U_a \}_{a \in 
\Phi}, T \bigr)$, and assume that the
Weyl group $W$ is infinite.
We denote by $G=\overline\Lambda^{{\hbox{\sevenrm Aut}}(\Delta)}$ the
corresponding topological Kac-Moody group defined in {\rm 1.B.2}, and 
by ${\Bbb A}$ the apartment of
$\Delta$ stabilized by $N$.
\rm\pec

R{\sevenrm EMARKS}. ---
1) By [R\'e1, 12.5.4], for each root $a$ the group $U_a$ is 
isomorphic to the ${\bf F}_q$-points of the
unipotent radical of a Borel subgroup in a rank one finite group of Lie type.
The finiteness assumption on the ground field implies the first 
assumption of 1.B.1,
and we have the lower bound $\mid \! U_a \! \mid \, \, \geq q$ for 
every root $a$.

2) We use the notation $T$ instead of $H$ for the twin root datum 
because $T$ is the group of rational
points of an ${\bf F}_q$-torus [R\'e1, 13.2.2].

\pec

E{\sevenrm XAMPLES}. --- 1) The group ${\rm SL}_2 \bigl( {\bf F}_q 
[t,t^{-1}] \bigr)$ is a split Kac-Moody group over
${\bf F}_q$ with infinite dihedral Weyl group.
The associated buildings are homogeneous trees of valencies
$1+q$ because the root groups of the affine twin Tits systems all 
have order $q$.

2) The group ${\rm SU}_3 \bigl( {\bf F}_q [t,t^{-1}] \bigr)$ is an 
almost split Kac-Moody group
over ${\bf F}_q$ with infinite dihedral Weyl group, too, but the 
associated buildings are semihomogeneous trees of valencies
$1+q$ and $1+q^3$ since the root groups have order $q$ or $q^3$ [R\'e4, 3.5].

\pec

In order to get lattices by 1.B.3, the condition requiring that 
$q_{\hbox{\sevenrm min}}$ be large enough, simply means that
$q$ is large enough.
In the case of split Kac-Moody groups, this assumption admits a sharp 
quantitative version
relating the growth series of the Weyl group and $q$.
Once the Haar measure of the full automorphism group of a building is 
normalized in such a
way that a chamber fixator has measure 1, the covolume of the Borel 
subgroup $\Gamma$ is
$\displaystyle \sum_{n \geq 0} {d_n \over q^n}=\sum_{w \in 
W}q^{-\ell(w)}$, where $d_n$ is the
number of elements in $W$ of length $n$ [R\'e2].

\pec

R{\sevenrm EMARKS}. --- 1) In the affine case, Kac-Moody groups are 
special cases of arithmetic groups over
function fields.
For this class of groups, a classical theorem by H. Behr and G. 
Harder implies that the
assumption on $q$ is empty [Mar, I.3.2.4].
Therefore an affine Kac-Moody group over a finite field is always a 
lattice of the product
of its two Euclidean buildings.

2) Let us consider a Kac-Moody group whose associated buildings have 
hyperbolic right-angled
regular $r$-gons as chambers -- see 2.E and [R\'e4, 4.1] for the 
existence of such groups.
Then the apartments are tilings of the hyperbolic plane ${\Bbb H}^2$ 
and the Weyl group $W$
has exponential growth.
More precisely, the growth series of $W$ is
$\displaystyle 1 + rt + \sum_{n \geq 2} r(r-2)^{n-1}t^n$, and the finiteness of
$\sum_n d_nq^{-n}$ amounts to $q \geq r-1$.
Therefore the Borel subgroup of such a Kac-Moody group is not always a lattice.

\pec

{\bf 1.C.2 } The inclusion $\Lambda<G$ implies that $G$ also admits a
$BN$-pair, since it is strongly transitive on the building $\Delta$ [Ro1 \S5].
We will prove a stronger result stressing the analogy between $G$ and 
semisimple groups over
local fields of positive characteristic.

\pec

T{\sevenrm HEOREM}.--- \it
Let $\Lambda$ and $G$ be groups as in {\rm 1.C.1}. \pec

{\rm (i) } The sixtuple $(G,N,U_c,\Gamma,T,S)$ defines a structure of 
refined Tits system for the associated topological
Kac-Moody group $G$.

{\rm (ii) } Let $F$ be a facet in the apartment ${\Bbb A}$.
Then the group $M({\Bbb A},F)$ of {\rm 1.A.4~} is finite of Lie type, 
$U_F$ is a pro-$p$ group, and the
parahoric subgroup $G_F = {\rm Fix}_{G}(F)$ decomposes as $G_F = 
M({\Bbb A},F) \ltimes U_F$.
\rm\pec

The proof will be given by a series of lemmas in 1.C.4 to 1.C.6.

\pec

Since the group $\Lambda$ acts by type-preserving isometries on 
$\Delta$, the stabilizer
of a facet $F$ in
$\Lambda$ is also its fixator.
We will denote it by $\Lambda(F)$; it is contained in $G_F$ but has 
no topological structure.
As explained in 1.A.4, any twin apartment ${\Bbb A}$ containing $F$ 
gives rise to a Levi
decomposition $\Lambda(F) = M({\Bbb A},F) \ltimes U(F)$.
Here the subgroup $M({\Bbb A},F)$ is abstractly isomorphic to the 
Kac-Moody group associated to a submatrix of the
generalized Cartan matrix defining $\Lambda$.
Since all facets in the metric buildings are spherical, the walls in 
${\Bbb A}$ containing
$F$ are finite in number; in fact, $M({\Bbb A},F)$ is a finite group 
of Lie type over ${\bf F}_q$.

\pec

R{\sevenrm EMARK}. --- It follows from theorem 1.C.2 that the group 
${\rm SL}_n\bigl( {\bf F}_q(\!(t)\!)
\bigr)$ satisfies the axioms of a refined Tits system, and that its 
parahoric subgroups admit semidirect product
decompositions. This splitting of facet fixators is not true for 
algebraic groups over local
fields of characteristic 0.
It is related to the existence of finite root groups with nice 
properties, which geometrically amounts to
the {\it Moufang property} -- see [Ro1, 6.4].
Other arguments (involving torsion, for instance) explaining why the 
analogy with the characteristic 0 case
is {\it not } relevant, will be given in 5.A and 5.B.

\pec

{\bf 1.C.3 } We choose a facet $F$ in ${\Bbb A}$, and define an 
exhaustion $\{ E_n \}_{n \geq 1}$
of $\Delta$, with respect to $F$.
We work inductively outwards from $F$.
The first term $E_1$ is (the closure of) ${\rm st}(F)$, the star of 
$F$ in $\Delta$.
To define further terms, suppose that $E_n$ is already defined.
Write ${\Bbb A}_n:= E_n \cap {\Bbb A}$ and choose a chamber $c_n$ of 
${\Bbb A}$ sharing a
panel with the boundary of ${\Bbb A}_n$.
Then we define ${\Bbb A}_{n+1}$ to be the convex hull of $c_n$ and 
${\Bbb A}_n$, and $E_{n+1}$ to be
$\Lambda(F).{\Bbb A}_{n+1}$, that is the set of all 
$\Lambda(F)$-transforms of ${\Bbb A}_{n+1}$.
By the Bruhat decomposition, for any chamber $c \! \in \! {\rm 
st}(F)$, ${\Bbb A}$ is a
complete set of representatives for the action of $U(c)<\Lambda(F)$ 
on $\Delta$.
Hence $\{ E_n \}_{n \geq 1}$ exhausts $\Delta$.
For each $n \geq 1$, the set $E_n$ is $\Lambda(F)$-stable.
We can write the groups $G_F$ and $U_F$ defined in 1.B.2 as projective limits:

\pec

\centerline{$G_F = \limproj_{n \geq 1} \Lambda(F) \! \mid_{E_n}$
\quad and \quad $U_F = \limproj_{n \geq 1} U(F) \! \mid_{E_n}$.}

\pec

We turn now to the proof of the theorem in 1.C.2.
We shall first prove, in 1.C.4, that $U_F$ is a pro-$p$ group, then 
in 1.C.5 we prove
the assertion about semi-direct products for parahoric subgroups, and 
finally in 1.C.6 we
deal with the refined Tits system for $G$.

\pec

{\bf 1.C.4 } Recall that each root group $U_a<\Lambda$ is a $p$-group (1.C.1).

\pec

L{\sevenrm EMMA}.--- \it
Let $\Lambda$ and $G$ be groups as in {\rm 1.C.1}. Then the group 
$U_F$ is pro-$p$.
\rm\pec

{\it Proof}. Recall that if a group $G$ acts on a set $S$, the 
notation $G \! \mid_S$ means
the group obtained from $G$ by factoring out the kernel of the action on $S$.
For each $n \geq 1$, the group $U(F) \! \mid_{E_{n+1}}$ when 
restricted to $E_n$ has a kernel
$K_n$.
This is represented by the following exact sequence:

\pec

\centerline{$1 \longrightarrow K_n \longrightarrow U(F) \! \mid_{E_{n+1}}
\longrightarrow U(F) \! \mid_{E_n} \longrightarrow 1$.}

\pec

The kernel $K_n$ is the fixator
${\rm Fix}_{U(F) \mid_{E_{n+1}}} (E_n)$ of $E_n$ in the restricted group
$U(F) \! \mid_{E_{n+1}}$ .
Since $U(F)$ fixes ${\rm st}(F)$, we have $U(F) \! \mid_{E_1} =
\{ 1 \}$.
By definition $U_F = \limproj_{n \geq 1} U(F) \! \mid_{E_n}$, so in 
order to show that $U_F$ is
pro-$p$, it is enough, by induction, to show that $K_n$ is a 
$p$-group for each $n \geq 1$.

\pec

Let us fix $n \geq 1$, and $u \! \in \! K_n$.
Let $\Pi_n$ be the panel of $c_n$ in the boundary of ${\Bbb A}_n$, 
and let $d_n$ be the other chamber
of ${\Bbb A}$ having $\Pi_n$ as a panel.
Then $d_n \! \in \! {\Bbb A}_n$, so $u$ fixes $d_n$, and by the 
Moufang property there
exists $v \! \in \! U_a$ such that $v^{-1}u.c_n=c_n$.
The element $v^{-1}u$ is in the \og unipotent radical\fg of the Borel subgroup
${\rm Fix}_\Lambda(d_n)$ and fixes $c_n \cup d_n$: hence it must fix 
the whole star of
$\Pi_n$.
Consequently, $u$ and $v$ coincide on ${\rm st}(\Pi_n)$ and thus
$u^{q_{\hbox{\sevenrm max}}}$ acts trivially on it.
In particular, $u^{q_{\hbox{\sevenrm max}}}$ fixes ${\Bbb A}_n \cup 
\{ c_n \}$ hence ${\Bbb
A}_{n+1}$ by convexity.

\pec

Every chamber of $E_{n+1} \setminus E_n$ is of the form $v.d$ for some
$d$ in the convex hull of ${\Bbb A}_n$ and $c_n$, and $v \! \in \! 
\Lambda(F)$. Then
$u^{q_{\hbox{\sevenrm max}}}(vd) = v(v^{-1}u^{q_{\hbox{\sevenrm max}}}v.d)
=v \bigl( (v^{-1}uv)^{q_{\hbox{\sevenrm max}}}.d \bigr)$.
Recall that $U(F)$ is normalized by $\Lambda(F)$, so that
applying the result of the above paragraph to $v^{-1}uv$ instead of
$u$, we get $(v^{-1}uv)^{q_{\hbox{\sevenrm max}}}.d=d$, and thus 
$u^{q_{\hbox{\sevenrm
max}}}(vd)=vd$.
This shows that the order of each element $u \! \in \! K_n$
divides $q_{\hbox{\sevenrm max}}$, so $K_n$ is a $p$-group.
Therefore $U_F$ is pro-$p$, as required.
\qed\pec

{\bf 1.C.5 } The following simple lemma will also be useful for the 
last part of the proof.

\pec

L{\sevenrm EMMA}.--- \it
Let $\Lambda$ and $G$ be groups as in {\rm 1.C.1}. \pec

{\rm (i) } We have $\Lambda\cap U_F = U(F)$, hence $M({\Bbb A},F) 
\cap U_F = \{ 1 \}$.
In particular, for any chamber $c$ in ${\Bbb A}$, we have $T \cap U_c 
= \{ 1 \}$.

{\rm (ii) } The group $G_F$ decomposes as $G_F = M({\Bbb A},F) \ltimes U_F$.
\rm\pec

{\it Proof}. (i).
First $\Lambda \cap U_F<{\rm Fix}_{\Lambda}(F) = \Lambda(F)$.
Thus if $g$ is in $\Lambda \cap U_F$, the Levi decomposition of $\Lambda(F)$
gives $g = mu$, with $m \! \in \! M({\Bbb A},F)$ and $u \! \in \! U(F)$.
Since $U_F$ fixes ${\rm st}(F)$, the set of chambers containing $F$, 
we see that
both $g$ and $u$, hence also $m$, fix ${\rm st}(F)$.
Therefore $m$ is central in the finite (reductive) group $M({\Bbb 
A},F)$ of Lie type; this implies
that $m$ lies in a torus, so its order is not divisible by $p$.
But $m$ is also a torsion element in the pro-$p$ group $U_F$, and is 
hence trivial.
The other assertions come from the trivial intersection $M({\Bbb 
A},F) \cap U(F)
= \{ 1 \}$ for any $F$.

\pec

(ii). For each $n \geq 1$, $E_n$ is $\Lambda(F)$-stable by definition, hence
$M({\Bbb A},F)$-stable. Since $M({\Bbb A},F)$ normalizes $U(F)$, it 
normalizes $U(F) \!
\mid_{E_n}$ in
$\Lambda(F) \! \mid_{E_n}$ for each $n \geq 1$. Hence $U_F$ is
normalized by $M({\Bbb A},F)$.

\pec

Let $g \! \in \! \Lambda(F)$. It can be written as $g = \limproj_{n \geq
1} g_n$,
with $g_n \! \in \! \Lambda(F) \! \mid_{E_n}$.
According to the Levi decomposition of $\Lambda(F)$, we have
$g_n=m_nu_n$ with
$m_n \! \in \! M({\Bbb A},F) \! \mid_{E_n}$ and $u_n \! \in \! U(F) \!
\mid_{E_n}$.
The group $M({\Bbb A},F)$ is finite, so up to extracting a 
subsequence, we may (and do)
assume that $m_n$ is a constant element $m$ in $M({\Bbb A},F)$.
This enables to write $g = m \cdot (\limproj_{n \geq 1} u_n)$, and
proves $G_F = M({\Bbb A},F) \cdot U_F$.

\pec

Finally (i) implies the trivial intersection $M({\Bbb A},F) \cap U_F 
= \{ 1 \}$, hence the
semidirect product decomposition $G_F = M({\Bbb A},F) \ltimes U_F$.
\qed\pec

{\bf 1.C.6 } We can now complete the proof of the theorem in 1.C.2.

\pec

L{\sevenrm EMMA}.--- \it
Let $\Lambda$ and $G$ be groups as in {\rm 1.C.1}. The sixtuple 
$(G,N,U_c,\Gamma,T,S)$ is a refined
Tits system as defined by the {\rm (RT)} axioms in {\rm 1.A.5}.
\rm\pec

{\it Proof}. We already know that $G$ admits a structure of $BN$-pair 
by strong transitivity of the $\Lambda$-action on
$\Delta$.
The main result we use for the verification below is [R\'e1, 1.5.4].
This theorem says that if
$\bigl( \Lambda, \{U_a \}_{a \in \Phi}, T \bigr)$ is a twin root 
datum indexed by the Coxeter system $(W,S)$, then
$(\Lambda,N,U_+,U_-,T,S)$ is a refined Tits system with Weyl group $W$.
Moreover we have $U_+ \cap s^{-1} U_- s=U_{a_s}$ for any $s \! \in \! S$.

\pec

Axiom (RT1). The group $T$ normalizes $U(c)$, and passing to the 
projective limit it
normalizes $U_c$. The other statements in axiom (RT1) are either 
clear or true from the
refined Tits system axioms for $\Lambda$.

\pec

Axioms (RT2). Define $U(s):= U_c \cap s^{-1}\Gamma s$.
The first point of the lemma in 1.C.5 implies $U(s)<U(c) \cap 
s^{-1}\Gamma s$, hence
$U(s)$ is the root group $U_{a_s}$
according to the relation of a refined Tits system and the twin root
datum of a Kac-Moody group [R\'e1, 1.6].
Axioms (RT2a) and (RT2b) are then clear because they just involve
subgroups of $\Lambda$.
For axiom (RT2c), we need $U_c = U_{a_s} \cdot (U_c \cap s^{-1}
U_c s)$. An element $u \! \in \! U_c$ fixes $c$, hence permutes the chambers
sharing their panel of type $s$ with
$c$. Since $U_{a_s}$ fixes $c$ and is (simply) transitive on the set of
chambers $\neq c$ of the
latter type, there exists an element $v \! \in \! U_{a_s}$ such that
$v^{-1}u$ fixes both $c$ and
$sc$. This proves (RT2c).

\pec

Axiom (RT3). Suppose now we have $\gamma n u_c = 1$, with $\gamma \! \in
\! \Gamma$, $n \! \in \!
N$ and $u_c \! \in \! U_c$. Then $\gamma n = u_c^{-1}$ belongs to $U_c
\cap \Lambda$, which
is $U(c)$ by the last lemma. Then all the factors are in $\Lambda$ and
we just have to apply
axiom (RT3) from the refined Tits system structure of the latter group.
\qed\pec

This completes the proof of the theorem in 1.C.2.

\pec

R{\sevenrm EMARKS}. --- 1) The groups $U_F$ are analogues of 
congruence subgroups in the classical case.

2) In the context of refined Tits systems, we have a refined Bruhat 
decomposition.
Hence, we could use the same arguments as in 1.B.3 to prove
\pec

L{\sevenrm EMMA}.--- \it
All parahoric subgroups of $G$ are commensurable.
\qed\rm\pec

Subsect. 1.C supports the analogy between topological Kac-Moody 
groups and semisimple
algebraic groups over non-Archimedean local fields of positive characteristic.
The analogy with the wider class of topological groups of Kac-Moody 
type is postponed to
the last (fifth) section.
In the next three sections, we develop the construction of twin root 
data exotic enough to
show that the analogy is strictly weaker in the wider context.

\vfill\eject

\centerline{\bf 2. Sketch of construction. Auxiliary computations}

\pec

We give here a summary of the construction, and introduce notation in
${\rm SL}_2$ in order to state relations to be used later.

\pec

{\bf 2.A } {\it Sketch of construction}. --- This construction will 
apply to trees as well as to some two-dimensional hyperbolic 
buildings.
We assume we are given the tiling of ${\bf R}$ by the segments
defined by the integers, or a
tiling of the hyperbolic plane ${\Bbb H}^2$ by regular right-angled
$r$-gons.
Here are the main steps.

\pec

1) Define the \og Borel subgroup\fg as the limit of an inductive system
described by the tiling.

2) To each type of vertex associate a \og unipotent radical\fg and a \og
Levi factor\fg.

3) Show that the \og Levi factor\fg admits an action on the \og
unipotent radical\fg.

4) In the two-dimensional case, define the other parabolic subgroups.

5) Amalgamate these groups along the inclusions given by the closure of
a chamber.

6) Verify the twin root datum axioms for the so-obtained amalgam $\Lambda$.

\pec

In order to apply this procedure, we will use relations in ${\rm 
SL}_2$ over a field.
The rest of the section introduces notation and relations in this
group.

\pec

R{\sevenrm EMARK}. --- Step 5) corresponds to defining $\Lambda$ as 
the fundamental group of a
complex of groups in the sense of H\ae fliger [BH, III.${\cal C}$] -- 
see 4.E.1 for more detailed explanations.

\pec

{\bf 2.B } {\it Notation in ${\rm SL}_2$}. --- We use Tits' 
convention for the parametrization of the negative root
group, and set:

\pec

\centerline{$u_i(r):= \pmatrix{1 & r \cr 0 & 1} \quad {\rm and } 
\quad u_{-i}(r):= \pmatrix{1 & 0 \cr -r & 1}$}

\pec

\centerline{$m_i(\lambda):= u_{-i}(\lambda^{-1})u_i(\lambda)u_{-i}(\lambda^{-
1})
= \pmatrix{0 & \lambda \cr -\lambda^{-1} & 0}$}

\pec

\centerline{$m_{-i}(\lambda):= u_i(\lambda^{-1})u_{-i}(\lambda)u_i(\lambda^{-
1})
= \pmatrix{0 & \lambda^{-1} \cr -\lambda & 0}$}

\pec

\centerline{$m_i:= m_i(1) \quad {\rm and } \quad t_i(\lambda):= 
m_i(\lambda)m_i^{-1}
= \pmatrix{\lambda & 0 \cr 0 & \lambda^{-1}}$}

\pec

{\bf 2.C } {\it Some straightforward relations}. --- The 
verifications are straightforward.

\pec

\centerline{$m_i(\lambda)^{-1} = m_i(-\lambda),
\quad m_{-i}(\lambda) = m_i(\lambda^{-1})
\quad {\rm and } \quad m_i^2 = -{\rm id}$;}

\pec

\centerline{$m_i(\lambda)u_i(r)m_i(\lambda)^{-1} = u_{-i}(\lambda^{-2}r)
\quad {\rm and }
\quad m_i(\lambda)u_{-i}(r)m_i(\lambda)^{-1} = u_i(\lambda^2 r)$;}

\pec

\centerline{$m_i(\lambda)t_i(\mu)m_i(\lambda)^{-1} = t_i(\mu^{-1}) = 
t_i(\mu)^{-
1}
\quad {\rm and } \quad m_i(\lambda)m_i(\mu) = t_i(-\lambda\mu^{-1})$;}

\pec

\centerline{$t_i(\lambda)u_i(r)t_i(\lambda)^{-1} = u_i(\lambda^2r)$.}

\pec

{\bf 2.D } {\it The product formula}. ---
We want to describe the product of two elements $g$ and $h$ in ${\rm 
SL}_2$ in terms of the Bruhat
decomposition.

\pec

C{\sevenrm ONVENTION}.--- \it
When $g$ is in the big cell ${\rm SL}_2 \pmatrix{* & * \cr \neq 0 & *}$,
we write it $g:= u_i(r)m_i(\lambda)u_i(r')$; when it is in the Borel
subgroup
${\rm SL}_2 \pmatrix{* & * \cr 0 & *}$, we write it $g:=
u_i(r)t_i(\lambda)$.
Similarly, when $h$ is in the big cell,
we write it $h:= u_i(s)m_i(\mu)u_i(s')$; when it is in the Borel
subgroup,
we write it $h:= u_i(s)t_i(\mu)$.
\rm\pec

{\bf 2.D.1 } First case: \og big cell $\cdot$ big cell $\in$ big
cell\fg, that is $r' + s \neq 0$.
Then:

$$gh = u_i \bigl( r-\lambda^2(r'+s)^{-1} \bigr) m_i \bigl( -\lambda\mu
(r'+s)^{-1} \bigr)
u_i \bigl( s'-\mu^2(r'+s)^{-1} \bigr).$$

\pec

{\bf 2.D.2 } Second case: \og big cell $\cdot$ big cell $\in$ Borel\fg,
that is $r' + s = 0$.
Then:

$$gh = u_i(r)t_i(-\lambda\mu^{-1})u_i(s') =
u_i(r+\lambda^2\mu^{-2}s')t_i(-\lambda\mu^{-1}).$$

\pec

{\bf 2.D.3 } Third case: \og Borel $\cdot$ big cell $\in$ big cell\fg.
Then: $gh = u_i(r+\lambda^2s)m_i(\lambda\mu)u_i(s')$.

\pec

{\bf 2.D.4 } Fourth case: \og big cell $\cdot$ Borel $\in$ big cell\fg.
Then: $gh = u_i(r)m_i(\lambda\mu^{-1})u_i(\mu^{-2}(r'+s) \bigr)$.

\pec

{\bf 2.D.5 } Fifth case: \og Borel $\cdot$ Borel $\in$ Borel\fg.
Then: $gh = u_i(r+\lambda^2s)t_i(\lambda\mu)$.

\pec

{\bf 2.E } {\it The examples from Kac-Moody theory}. ---
In this subsection, we recall some basic facts on Kac-Moody groups, 
and we present a very specific
case of Kac-Moody groups which admit a well-defined ground field.
Our aim, later in sect. 4, is to construct groups of the same type 
(i.e. having the same Weyl group) but
using several different ground fields.

\pec

{\bf 2.E.1} A {\it generalized Cartan matrix } is an integral matrix 
$A=[A_{ij}]_{i, j \in I}$ such that
$A_{ii}=2$, $A_{ij} \leq 0$ whenever $i \neq j$ and $A_{ij}=0 
\Leftrightarrow A_{ji}=0$.
A {\it Kac-Moody root datum } with generalized Cartan matrix $A$ is a $5$-tuple
${\cal D}=\bigl(I, A, X^*, \{a_i\}_{s \in I}, \{h_i\}_{s \in I} 
\bigr)$, where $X^*$ and $X_*$ are free
${\bf Z}$-modules ${\bf Z}$-dual to one another and such that the 
elements $a_i$ of $X^*$ and $h_j$
of $X_*$ satisfy $a_i(h_j)=A_{ij}$ for all $i, j \! \in \! I$ [T2, 
introduction].
The datum ${\cal D}$ is called {\it simply connected } if the 
elements $h_i$ freely generate $X_*$.
A Kac-Moody root datum is a short way to encode the defining 
relations of a Kac-Moody group functor
as defined by J. Tits [T2, 3.6].
When $A$ is a (classical) Cartan matrix, such a functor is the 
functor of points of a split reductive
group scheme over ${\bf Z}$, and $X^*$ (resp. $X_*$) is the group of 
characters (resp. cocharacters)
of a maximal split torus.

\pec

Let us denote by $\Lambda$ the value over a field ${\bf K}$ of
the Kac-Moody group functor associated to a Kac-Moody root datum ${\cal D}$.
We mentioned in the introduction of 1.C that $\Lambda$ admits a twin root datum
$\bigl( \Lambda, \{ U_a \}_{a \in \Phi}, T \bigr)$.
We describe it more precisely now, and refer to [R\'e1 \S 8] for 
proofs or references.
The group $T$ is the split torus $T={\rm Hom}({\bf K}[X^*],{\bf K})$, 
with cocharacter group $X_*$.
For each $i \! \in \! I$, there is an involution $s_i$ of the ${\bf Z}$-module
$\displaystyle \bigoplus_{i \in I} {\bf Z}a_i$ defined by $s_i(a_j) = 
a_j - A_{ij}a_i$.
The group generated by these involutions is a Coxeter group $W$, 
called the {\it Weyl group }of
$\Lambda$.
In the Kac-Moody case, the root system $\Phi$ of $W$, as defined in 
1.A.1, is in bijection with the
elements of the form $wa_i$.
The roots $a \! \in \! \Phi$ index a family of isomorphisms
$u_a : ({\bf K},+) \simeq U_a$, where $U_a$ is a root group of $\Lambda$.
In the case of a simply connected root datum ${\cal D}$, the
standard torus $T$ is the product $\displaystyle \prod_{i \in I} 
h_i({\bf K}^\times)$ of
multiplicative one-parameter subgroups.
The coordinates of a root $a$ give the powers of the multiplicative 
parameters by which $T$ operates
on the root group $U_a$.
Writing $t=\prod_{i \in I}h_i(\lambda_i) \! \in \! T$, we have:

\pec

\centerline{$t \cdot u_a(r) \cdot t^{-1} = u_a \bigl( \prod_{i \in I} 
\lambda_i^{a(h_i)} r)$.}

\pec

The commutation relations between root groups can be made quite 
explicit by computations \og \`a la
Steinberg\fg in a ${\bf Z}$-form of the universal enveloping algebra 
of the Lie algebra attached to
${\cal D}$.
This is done for instance in [R\'e1 \S 9] in order to define an 
adjoint representation for Kac-Moody
groups.

\pec

{\bf 2.E.2} We now turn to the special cases we are interested in.
Pick a prime power $q \geq 3$ and consider the generalized Cartan 
matrix $A$ indexed by
${\bf Z}/r$ ($r \geq 5$) in which $A_{i, i+1}=0$ and $A_{i j}=1-q$
for $j \neq i, i\pm1$.
Then the associated Weyl group is the hyperbolic reflection group 
arising from the tiling of the
hyperbolic plane by regular right-angled $r$-gons [R\'e1 \S 13].
The root system of the tiling is studied geometrically in 4.C, where 
it is shown in particular that roots
have a well-defined type: for instance, $wa_{i}$ has type $i$.

\pec

We choose the ground field to be ${\bf F}_q$, in which case the 
choice of off-diagonal coefficients makes the action of
tori on root groups as trivial as possible (using the fact that the 
$(q-1)$-th power of any element in
${\bf F}_q^\times$ is 1).
By induction each root $wa_i$ can be written $\pm a_{i} + (q-1) \, 
\times$ a linear
combination of simple roots, and a multiple of $q-1$ leads to a 
trivial action. Therefore the action of the
standard torus $T$ on the root group $U_a$ is simply multiplication
by the square (or the inverse of the square) of the multiplicative 
parameter having the same type as
$a$.

\pec

We now turn to commutation relations between root groups.
Since $q \geq 3$, J.-Y. H\'ee's work [Cho, 5.8] shows that
a pair of roots where one contains the other leads to a trivial 
commutation relation.
In the tiling of ${\Bbb H}^2$ any two walls are parallel or orthogonal.
When two roots have orthogonal walls, the corresponding root groups 
commute, so the only
possibility to have two non-commuting root groups is when the two 
roots (or their two
opposites) intersect along a strip in ${\Bbb H}^2$.
In that case, there is no relation at all, because the pair is not 
prenilpotent [T4]: the
free product ${\bf F}_q * {\bf F}_q$ of the root groups injects in 
the Kac-Moody group
$\Lambda$.

\pec

R{\sevenrm EMARK}. ---
Our coefficients in the generalized Cartan matrix were suitably chosen with
respect to the characteristic of the ground-field.
This trick is also the starting point of Ree-Suzuki torsions, which 
were defined for Kac-Moody groups by J.-Y. H\'ee [H].

\pec

The above special cases of Kac-Moody groups show the existence of 
automorphism groups
of Moufang twinnings with particularly simple commutation relations 
(between root groups) and
actions of tori on root groups.
In these cases, root groups and multiplicative one-parameter 
subgroups have a well-defined
type, and when the types are different, commutation relations are trivial.
This is our model for the constructions in sections 3 and 4.

\vfill\eject

\centerline{\bf 3. Twinning trees}

\pec

We now apply procedure 2.A to trees.
The computations for step 3) will be done in detail, enabling us to 
concentrate on more
conceptual arguments in the two-dimensional case of Fuchsian 
buildings -- see Sect. 4.

\pec

{\bf 3.A } {\it Geometric description of the root system}. ---
In this section, the Coxeter complex we are interested in is the 
tiling of the real line
${\bf R}$ by the length 1 segments whose vertices are integers.
We denote by $E$ the segment $[0,1]$, and by $s_0$ (resp. $s_1$) the 
reflection with
respect to $0$ (resp. $1$).
This provides a geometric realization of the Coxeter complex of the 
infinite dihedral group
$D_\infty$, seen as the group generated by $s_0$ and $s_1$.
The roots are the half-lines defined by the integers, the positive 
roots are by definition
those which contain $E$.
We denote by $a_0$ (resp. $a_1$) the positive root having boundary 
$0$ (resp. $1$).

\pec

P{\sevenrm ICTURE}. ---

\pec

\centerline{\epsfysize=20mm $$\epsfbox{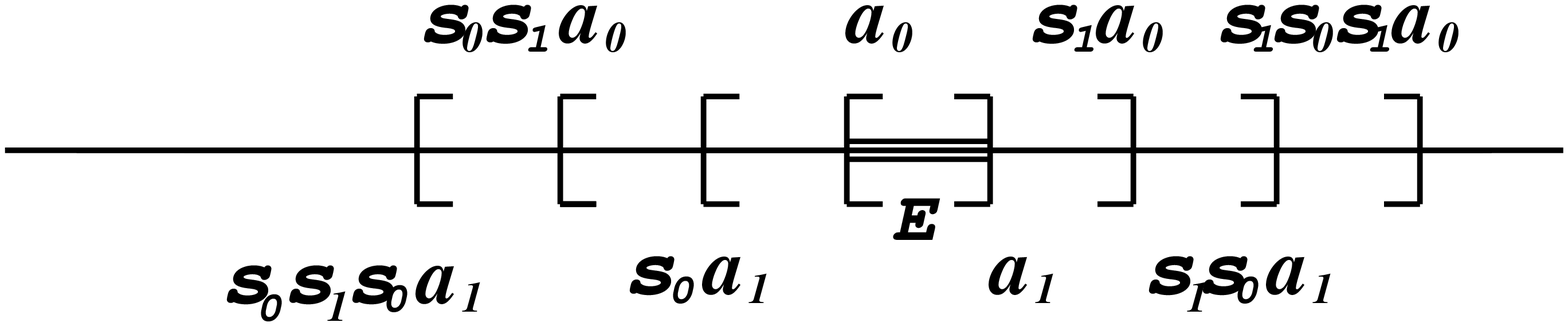}$$}

\pec

Each vertex has type 0 or 1, with the obvious notation, and the 
boundary of a root is called its
{\it vertex}.
The {\it type } of a root $a$ is the type of its vertex, and is 
denoted by $\iota(a)$.

\pec

The notion of a prenilpotent pair of roots in a Coxeter complex is 
defined in 1.A.1, and will be used later.
In this section, with the Coxeter complex as defined above, a 
positive root $a$ is prenilpotent with $a_i$ if
and only if it contains $a_i$.
Each positive root contains either $a_0$ or $a_1$, but not both, so 
the positive roots
fall into four disjoint subsets ${\cal P}(a_i,j)$ for $i,j \! \in \! 
\{ 0;1 \}$, where ${\cal P}(a_i,j)$ is the set
of positive roots of type $j$ which contain $a_i$.
Note that these subsets of positive roots can also be defined as follows:

\pec

\centerline{${\cal P}(a_0,0):= \{ (s_0 s_1)^m a_0 \}_{m \geq 0},
\quad {\cal P}(a_0,1):= \{ (s_0 s_1)^m s_0 a_1 \}_{m \geq 0}$,}

\pec

\centerline{${\cal P}(a_1,0):= \{ (s_1 s_0)^m s_1 a_0 \}_{m \geq 0},
\quad {\cal P}(a_1,1):= \{ (s_1s_0)^m a_1 \}_{m \geq 0}$.}

\pec

We shall use ${\cal P}(a_i)$ to mean ${\cal P}(a_i,0)\cup{\cal 
P}(a_i,1)$, in other words the set of
positive roots containing $a_i$ (i.e. prenilpotent with $a_i$).

\pec

{\bf 3.B } {\it The Borel subgroup. Unipotent subgroups}. --- We pick 
two fields ${\bf K}_0$ and ${\bf K}_1$.
To each positive root $a$ of type $i$ is attached a copy of the additive
group $({\bf K}_i,+)$.
We denote it by $U_a:= \{ u_a(k): k \! \in \! {\bf K}_i \}$
($u_a$ is the chosen isomorphism between the root group and its field).

\pec

D{\sevenrm EFINITION}.--- \it
{\rm (i) } Let the abelian group $A_i$ be the direct sum of the root groups
associated to the roots in ${\cal P}(a_i)$.

{\rm (ii) } Let the group $U_+$ be the free product of $A_0$ and $A_1$, that is
$U_+:= A_0* A_1$.

{\rm (iii) } Let the group $T$ be the product ${\bf K}_0^\times 
\times {\bf K}_1^\times$ of the
multiplicative groups of the chosen fields.
\rm\pec

For a (positive) root $a$ of type $i$, we define $\epsilon_a \! \in 
\! \{ \pm 1 \}$ to
be $(-1)^m$, where $m$ is the exponent appearing in 3.A.
It is the parity of the number of vertices of type $i$ in the interior of the
segment joining the middle of $E$ and the vertex $\partial a$.
Using the notation of Sect. 2, $T$ is the maximal torus
$\{ t_i(\lambda) t_{1-i}(\mu): \lambda, \mu$ invertible$\}$ of ${\rm 
SL}_2({\bf K}_0) \times {\rm
SL}_2({\bf K}_1)$, and we define the action of $T$ on $U_+$ by:

\pec

\hskip 3mm (3B1)
\hskip 3mm $t_j(\lambda)u_a(k)t_j(\lambda)^{-1}:= 
u_a(\lambda^{2\epsilon_a\delta_{j\iota(a)}}k)$,

\pec

where the Kronecker symbol $\delta_{j\iota(a)}$ means that the 
element $t_j(\lambda)$ induces
multiplication of $k$ by $\lambda^{2\epsilon_a}$ only if the type of 
the root $a$ is $j$.
In particular, $t_j(\lambda)$ centralizes each root group of type $\neq j$.

\pec

D{\sevenrm EFINITION}.--- \it
The {\rm (positive) Borel subgroup $B$ } is the semi-direct product
$B:= T \ltimes U_+$.
\rm\pec

Let $V^i$ be the subgroup of $U_+$ generated by all the positive root
groups except $U_i$; it is not normal
in $U_+$.

\pec

D{\sevenrm EFINITION}.--- \it
The group $U^i$ is the normal closure of $V^i$ in $U_+$.
\rm\pec

In view of the free product structure of $U_+$, $U^i$ is generated by
$u_a(k)$ for
$a \! \in \! {\cal P}(a_i) \setminus \{ a_i \}$
and $u_i(r)u_a(k)u_i(r)^{-1}$ for $a \! \in \! {\cal P}(a_{1-i})$, with
$r \! \in \! {\bf K}_i$.

\pec

{\bf 3.C } {\it Actions of Levi factors on unipotent radicals}. ---
Let us start with the definition of another subgroup.

\pec

D{\sevenrm EFINITION}.--- \it
The {\rm (standard) Levi factor of type $i$} is the direct product

\pec

\centerline{$L_i:= {\rm SL}_2({\bf K}_i) \times {\bf K}_{1-i}^\times$.}
\rm\pec

In order to define a parabolic subgroup as a semi-direct product $L_i
\ltimes U^i$, we must
define an action of $L_i$ on the group $U^i$.
We define actions of ${\rm SL}_2({\bf K}_i)$ and ${\bf 
K}_{1-i}^\times$ which are easily seen to commute with one
another.
Hence we can deal with each factor separately.

\pec

The action of the torus ${\bf K}_{1-i}^\times$ on $U^i$ is that 
obtained as a subgroup of $T$.
This makes sense because $T$ obviously stabilizes $U^i$.
We turn now to the action of the factor ${\rm SL}_2({\bf K}_i)$.
In this subsection, we define the actions of the generators $u_i(r)$ 
and $m_i(\lambda)$
given in 2.B.

\pec

For the root groups, we set for $k$ in ${\bf K}_0$ or ${\bf K}_1$, and
$r,s$ in ${\bf K}_i$:

\pec

\hskip 3mm (3C1)
\hskip 3mm $u_i(r)u_a(k)u_i(r)^{-1}:= u_a(k)$
\hskip 3mm for $a \! \in \! {\cal P}(a_i) \setminus \{ a_i \}$,

\pec

\hskip 3mm (3C2)
\hskip 3mm $u_i(s) \bigl( u_i(r)u_a(k)u_i(r)^{-1} \bigr) 
u_i(s)^{-1}:=u_i(r+s)u_a(k)u_i(r+s)^{-1}$
\hskip 3mm for $a \! \in \! {\cal P}(a_{1-i})$.

\pec

The elements $m_i(\lambda)$, lifting the Weyl group reflections (see 
Sect. 2), have a conjugation action on the
positive root group elements $u_a(k)$, defined by:

\pec

\hskip 3mm (3C3)
\hskip 3mm$m_i(\lambda)u_a(k)m_i(\lambda)^{-1}
:=u_{s_ia}(\lambda^{-2\epsilon_a\delta_{i\iota(a)}}k)$.

\pec

We must also define a conjugation action of $m_i(\lambda)$ on 
elements $u_i(r)u_a(k)u_i(r)^{-1}$ in the free product
of $U_a$ and $U_i$, whenever $a$ is a positive root containing 
$a_{1-i}$, and $u_i(r)\neq 1$.
We set:

\pec

\hskip 3mm (3C4)
\hskip 3mm $\displaystyle
m_i(\lambda) \bigl( u_i(r)u_a(k)u_i(r)^{-1} \bigr) m_i(\lambda)^{-1}
:= u_i \bigl({-\lambda^{2} \over r}\bigr)
u_a\bigl( ({-\lambda \over r})^{2\epsilon_a\delta_{i\iota(a)}}k \bigr)
u_i \bigl({-\lambda^{2} \over r}\bigr)^{-1}$.

\pec

R{\sevenrm EMARK}. ---
The Kronecker symbol $\delta_{i\iota(a)}$ in the exponents 
${}^{2\epsilon_a\delta_{i\iota(a)}}$ involves
the types of roots $i$ and $\iota(a)$.
It simply means that the element $m_i(\lambda)$ induces a multiplication of the
additive parameter $k$ in $u_a(k)$ by a factor $\lambda^{-2\epsilon_a}$ or
$\displaystyle ({-\lambda \over r})^{2\epsilon_a}$ only if the type 
of the root $a$ is $i$.

\pec

{\bf 3.D } {\it Checking the product relation}. ---
In this subsection, we make sure that the individual actions
above define an action of $L_i$ on $U^i$.
We must show that given any two elements
$g, h$ of ${\rm SL}_2({\bf K}_i) \times {\bf K}_{1-i}$ and a
generator $v$ of $U^i$, we always have

\pec

\centerline{$g(hvh^{-1})g^{-1}=(gh)v(gh)^{-1}$.}

\pec

This equality involves the product formula for $g\cdot h$, so we must 
take into account the five cases
of 2.D.
The form of the generator $v$ of $U^i$ will also play a role.
A generator $v=u_i(t)u_a(k)u_i(t)^{-1}$ with $a$ a root containing 
$a_{1-i}$ (hence not prenilpotent
with $a_i$), will be referred to as a generator {\it of the first type}.
A generator $v=u_a(k)$, with $a$ a root containing $a_i$ (hence 
prenilpotent with $a_i$), will be referred
to as a generator {\it of the second type}.

\pec

{\bf 3.D.1 } \og big cell $\cdot$ big cell $\in$ big cell\fg-- see 
2.D.1, and recall that

\pec

\centerline{$g=u_i(r)m_i(\lambda)u_i(r')$ \qquad and \qquad 
$h=u_i(s)m_i(\mu)u_i(s')$.}

\pec

We introduce the notation $R=r'+s$ and $S=s'+t$ for the remainder of 
section 3.D.
In 3.D.1, as in 2.D.1, we have $R \neq 0$ and

\pec

\centerline{$\displaystyle
gh = u_i \bigl( r-{\lambda^2\over R} \bigr) m_i \bigl( {-\lambda\mu 
\over R} \bigr) u_i \bigl( s'-{\mu^2\over R} \bigr)$.}

\pec

Let us deal with a generator $v$ of the first type.
By (3C2) $hvh^{-1} = [u_i(s)m_i(\mu)u_i(S)]u_a(k)[\cdots]^{-1}$, and 
assuming that $S \neq 0$, we have by (3C4) and
(3C2):

\pec

\centerline{$\matrix{hvh^{-1}
&= &
u_i(s - \mu^2S^{-1}) u_a \bigr(( -\mu 
S^{-1})^{2\epsilon_a\delta_{i\iota(a)}}k \bigr)u_i(\cdots)^{-1}, 
\hfill \cr
g(hvh^{-1})g^{-1}
&= &
[u_i(r)m_i(\lambda)u_i(R-\mu^2 S^{-1})]u_a \bigr((
\mu^2S^{-1})^{2\epsilon_a\delta_{i\iota(a)}}k\bigr)[\cdots]^{-1}.\hfill}$}

\pec

Under the further assumption that $R-\mu^2 S^{-1} \neq 0$, (3C4) and 
(3C2) give:

\pec

\centerline{$\matrix{\displaystyle g(hvh^{-1})g^{-1}
&= & \displaystyle
u_i\bigl(r-{\lambda^2 \over R-\mu^2 S^{-1}}\bigr)\cdot
u_a\bigl(({-\lambda \over R-\mu^2 
S^{-1}})^{2\epsilon_a\delta_{i\iota(a)}}(-\mu 
S^{-1})^{2\epsilon_a\delta_{i\iota(a)}} k\bigr)\cdot
u_i\bigl(\cdots\bigr)^{-1} \hfill\cr
&= & \displaystyle
u_i\bigl(r-{\lambda^2S \over RS-\mu^2}\bigr) \cdot
u_a \bigl( ({\lambda\mu \over 
RS-\mu^2})^{2\epsilon_a\delta_{i\iota(a)}}k\bigr) \cdot 
u_i\bigl(\cdots\bigr)^{-1}. \hfill
}$}

\pec

Now assuming $(R-\mu^2 S^{-1}) \cdot S \neq 0$, we use (3C4) and 
(3C2) to conjugate $v$ by $gh$:

\pec

\centerline{$\matrix{(gh)v(gh)^{-1}
&= & \displaystyle
[u_i( r-\lambda^2R^{-1})\displaystyle m_i(-\lambda\mu R^{-1})u_i( S 
-\mu^2R^{-1})]u_a(k)[\cdots]^{-1}
\hfill \cr
&= & \displaystyle
u_i(r-{\lambda^2 \over R} - {\lambda^2\mu^2 \over R^2(S-\mu^2R^{-1})})\cdot
u_a\bigl((-{-\lambda\mu \over 
R(S-\mu^2R^{-1})}\bigr)^{2\epsilon_a\delta_{i\iota(a)}}k\bigr)\cdot
u_i(\cdots)^{-1} \hfill \cr
&= & \displaystyle
u_i\bigl(r-{\lambda^2S \over RS-\mu^2}\bigr) \cdot
u_a \bigl( ({\lambda\mu \over 
RS-\mu^2})^{2\epsilon_a\delta_{i\iota(a)}}k\bigr) \cdot
u_i\bigl(\cdots\bigr)^{-1},
\hfill}$}

\pec

which proves $g(hvh^{-1})g^{-1}=(gh)v(gh)^{-1}$ when $S\neq 0$ and 
$R-\mu^2 S^{-1} \neq 0$.

\pec

Now, suppose $S \neq 0$ but $R-\mu^2 S^{-1}=0$.
Then the first equation above for $g(hvh^{-1})g^{-1}$ simplifies, and 
using (3C3) and (3C1)
(which implies that $u_i$ commutes with $u_{s_ia}$) we obtain:

\pec

\centerline{$\displaystyle g(hvh^{-1})g^{-1} =
u_i(r) u_{s_ia}\bigl( ({\mu\over\lambda 
S})^{2\epsilon_a\delta_{i\iota(a)}} k \bigr) u_i(r)^{-1} =
\displaystyle u_{s_ia}\bigl( ({\mu\over\lambda 
S})^{2\epsilon_a\delta_{i\iota(a)}} k \bigr),$}

\pec

and by (3C3) and (3C1) again:

\pec

\centerline{$\displaystyle (gh)v(gh)^{-1}
= [u_i \bigl( r-\lambda^2 R^{-1} \bigr) m_i \bigl(-\lambda \mu R^{-1} 
\bigr)] u_a(k) [\cdots]^{-1}
= u_{s_ia}\bigl( ({R\over\lambda\mu})^{2\epsilon_a\delta_{i\iota(a)}} 
k \bigr)$.}

\pec

Using $\displaystyle R={\mu^2\over S}$, we obtain the desired equality
$g(hvh^{-1})g^{-1}=(gh)v(gh)^{-1}$ in this case.

\pec

Finally if $S = 0$, by (3C3) and (3C1) we have:

\pec

\centerline{$\displaystyle hvh^{-1}
= [u_i (s) m_i(\mu)] u_a(k) [\cdots]^{-1} = u_i(s)
u_{s_ia}({k\over \mu^{2\epsilon_a\delta_{i\iota(a)}}}) u_i(s)^{-1}
= u_{s_ia}({k\over \mu^{2\epsilon_a\delta_{i\iota(a)}}})$.}

\pec

By (3C3), 
$m_i(\lambda)u_{s_ia}(k)m_i(\lambda)^{-1}=u_a(\lambda^{2\epsilon_a\delta_{i\iota(a)}})$. 

Using this and (3C2):

\pec

\centerline{$\displaystyle\quad g(hvh^{-1})g^{-1} =
[u_i(r)m_i(\lambda)] u_{s_ia}(\mu^{-2\epsilon_a\delta_{i\iota(a)}} k) 
[\cdots]^{-1} =
u_i(r) \cdot
u_a \bigl( ({\lambda\over \mu})^{2\epsilon_a\delta_{i\iota(a)}} k \bigr)
\cdot u_i(r)^{-1}$.}

\pec

Again with $S=0$, conjugation of $v$ by $gh$ gives by (3C4) and (3C2):

\pec

\centerline{$\matrix{(gh)v(gh)^{-1}
&= & [u_i \bigl( r-\lambda^2 R^{-1} \bigr) m_i \bigl(-\lambda \mu
R^{-1} \bigr)
u_i \bigl( -\mu^2 R^{-1} \bigr)] u_a(k) [\cdots]^{-1} \hfill \cr
&= & \displaystyle [u_i \bigl( r- {\lambda^2 \over R} \bigr)
u_i \bigl( {-\lambda^2\mu^2 R^{-2} \over -\mu^2 R^{-1}}
\bigr)]
u_a \bigl( ({\lambda \mu R^{-1} \over -\mu^2 R^{-1}}
)^{2\epsilon_a\delta_{i\iota(a)}} k\bigr) [\cdots]^{-1} \hfill \cr
&= & \displaystyle u_i(r) \cdot u_a \bigl( ({\lambda\over 
\mu})^{2\epsilon_a\delta_{i\iota(a)}} k \bigr)
\cdot u_i(r)^{-1},\hfill}$}

\pec

hence $g(hvh^{-1})g^{-1}=(gh)v(gh)^{-1}$.

\pec

The case of a generator $v=u_a(k)$ of the second type is simpler, 
because $a \! \in \! {\cal P}_i-\{ a_i\}$, so by (3C1)
$u_i$ commutes with $u_a$.
Using (3C2), (3C3) and (3C4) one obtains:

\pec

\centerline{
$\displaystyle g(hvh^{-1})g^{-1}=(gh)v(gh)^{-1}=
u_i(r-{\lambda^2\over R})\cdot
u_{s_ia}\bigl(({R\over\lambda\mu})^{2\epsilon_a\delta_{i\iota(a)}}\bigr)\cdot
u_i(\cdots)^{-1}$.}

\pec

{\bf 3.D.2 } \og big cell $\cdot$ big cell $\in$ Borel\fg -- see 
2.D.2: with $g$ and $h$ as above, we now have $R=0$,
and $\displaystyle gh = u_i(r+{\lambda^2\over 
\mu^2}s')t_i({-\lambda\over \mu})$.
Let us consider a generator $v$ of the first type.
For $g(hvh^{-1})g^{-1}$, we compute as in 3.D.1 (with $R=0$), the 
case $S \neq 0$ but $R-\mu^2
S^{-1}=0$ being excluded.
We obtain in any case:
$\displaystyle g(hvh^{-1})g^{-1} =
u_i(r+{S\lambda^2\over\mu^2})\cdot
u_a \bigl( ({\lambda\over \mu})^{2\epsilon_a\delta_{i\iota(a)}} k \bigr)
\cdot u_i(\cdots \bigr)^{-1}$,
which equals $ (gh)v(gh)^{-1}$ by (3C2) and (3B1).

\pec

For a generator of the second type, we have:
$\displaystyle g(hvh^{-1})g^{-1} = (gh)v(gh)^{-1} =
u_a \bigl( ({\lambda\over\mu})^{2\epsilon_a\delta_{i\iota(a)}} k \bigr)$.

\pec

{\bf 3.D.3 } \og Borel $\cdot$ big cell $\in$ big cell\fg -- see 
2.D.3, and recall that

\pec

\centerline{$g=u_i(r)t_i(\lambda)$, \qquad $h=u_i(s)m_i(\mu)u_i(s')$ 
\qquad and \qquad
$gh = u_i(r+\lambda^2s)m_i(\lambda\mu)u_i(s')$.}

\pec

For a generator of the first type, we have for $S \neq 0$:

\pec

\centerline{
$\displaystyle g(hvh^{-1})g^{-1} = (gh)v(gh)^{-1} =
u_i(r+\lambda^2s-{\lambda^2\mu^2\over S})
u_a \bigl(({\lambda\mu\over S})^{2\epsilon_a\delta_{i\iota(a)}}k\bigr)
u_i(\cdots)^{-1}$,}

\pec

and for $S=0$: \qquad $\displaystyle g(hvh^{-1})g^{-1} = (gh)v(gh)^{-1} =
u_{s_ia}\bigl(({1\over\mu\lambda})^{2\epsilon_a\delta_{i\iota(a)}}k\bigr)$.

\pec

For a generator of the second type, we have:

\pec

\centerline{$\displaystyle g(hvh^{-1})g^{-1}=(gh)v(gh)^{-1}
=u_i(r + \lambda^2 s)
u_{s_ia}\bigl(({1\over\mu\lambda})^{2\epsilon_a\delta_{i\iota(a)}}k\bigr)
u_i(\cdots)^{-1}$.}

\pec

{\bf 3.D.4 } \og big cell $\cdot$ Borel $\in$ big cell\fg -- see 
2.D.4, and recall that

\pec

\centerline{$g=u_i(r)m_i(\lambda)u_i(r')$, \qquad $h=u_i(s)t_i(\mu)$ 
\qquad and \qquad
$\displaystyle gh = u_i(r)m_i({\lambda\over\mu})u_i({R\over \mu^2} \bigr)$.}

\pec

Write $T=R+\mu^2t$, and note that $R\neq 0$ in this case.

\pec

Let us deal with a generator of the first type.
Using (2C) and (3B1) we have:

\pec

\centerline{
$hvh^{-1}=u_i(s+\mu^2t) \cdot u_a(\mu^{2\epsilon_a\delta_{i\iota(a)}} 
k) \cdot u_i(\cdots)^{-1}$.}

\pec

Hence under the assumption that $T\neq 0$, we get by (3C4) and (3C2):

\pec

\centerline{
$\displaystyle g(hvh^{-1})g^{-1}
=u_i(r-{\lambda^2 \over T})\cdot
u_a\bigl(( {\lambda \mu \over 
T})^{2\epsilon_a\delta_{i\iota(a)}}k\bigr)\cdot u_i(\cdots)^{-1}$,}

\pec

which equals $(gh)v(gh)^{-1}$, since by (3C4) and (3C2) we have:

\pec

\centerline{$\displaystyle (gh)v(gh)^{-1}
= [u_i(r) m_i(\lambda\mu^{-1}) u_i \bigl({T \over \mu^2})]u_a(k) [\cdots]^{-1}
= [u_i(r) u_i({-\lambda^2 \over T})]u_a \Bigl(\bigl({\lambda \mu\over 
T}\bigr)^{2\epsilon_a\delta_{i\iota(a)}} k \Bigr)
[\cdots]^{-1}.$}

\pec

When $T=0$, we simply have:
$\displaystyle g(hvh^{-1})g^{-1}=(gh)v(gh)^{-1} =
u_{s_ia} \bigl( ({\mu \over \lambda})^{2\epsilon_a\delta_{i\iota(a)}}k \bigr)$.

\pec

For a generator of the second type, we have in any case:

\pec

\centerline{$\displaystyle g(hvh^{-1})g^{-1}=(gh)v(gh)^{-1}=
u_i(r)\cdot u_{s_ia} \bigl( 
({\mu\over\lambda})^{2\epsilon_a\delta_{i\iota(a)}} k\bigr)\cdot
u_i(\cdots)^{-1}$.}

\pec

{\bf 3.D.5 } \og Borel $\cdot$ Borel $\in$ Borel\fg -- see 2.D.5, and 
recall that
$g=u_i(r)t_i(\lambda)$, $h=u_i(s)t_i(\mu)$ and $gh = 
u_i(r+\lambda^2s)t_i(\lambda\mu)$.
The common value of $g(hvh^{-1})g^{-1}$ and of $(gh)v(gh)^{-1}$ is:

\pec

$\matrix{\hskip 6mm
u_i(r+\lambda^2s+(\lambda\mu)^2t) \cdot
u_a\bigl((\lambda\mu)^{2\epsilon_a\delta_{i\iota(a)}}k\bigr) \cdot 
u_i(\cdots)^{-1} \hfill
&\hbox{\rm for a generator of the first type,} \hfill\cr
\hskip 6mm
u_a\bigl((\lambda\mu)^{2\epsilon_a\delta_{i\iota(a)}}k\bigr) \hfill
&\hbox{\rm for a generator of the second type.}}$

\pec

This concludes the computations and enables us to introduce the following

\pec

D{\sevenrm EFINITION}.--- \it
The {\rm (standard) parabolic subgroup of type $i$ } is the semidirect product

\pec

\centerline{$P_i:= L_i \ltimes U^i$.}
\rm\pec

We can now turn to combinatorial considerations.

\pec

{\bf 3.E } {\it Group combinatorics}. --- This subsection is 
dedicated to our main constructive result about non Kac-Moody Moufang 
twin trees.

\pec

{\bf 3.E.1 } We can now define the group we are interested in.

\pec

D{\sevenrm EFINITION}.--- \it
The {\rm (abstract) group of Kac-Moody type } associated to the 
choices of fields
above is the amalgam $\Lambda:= P_0 *_{B} P_1$.
\rm\pec

This amalgam makes sense because we have $P_i = L_i \ltimes U^i$,
with $L_i = {\rm SL}_2({\bf K}_i) \times {\bf K}_{1-i}^\times$.
In view of the actions of the subgroups $T:= {\bf K}_{i}^\times 
\times {\bf K}_{1-i}^\times$ and
$U_i$ of $L_i$ on $U^i$, the group $TU_i \ltimes U^i$ is isomorphic to
(and identified with) $B$.
By the Bruhat decomposition for ${\rm SL}_2$, each parabolic subgroup 
admits the following decomposition:

\pec

\centerline{$(\star) \quad P_i
= B \sqcup U_i m_i B \quad (i \! \in \! \{ 0;1 \})$.}

\pec

Consequently, an element $\gamma$ of $\Lambda$ can be written

\pec

\centerline{$\gamma = u_1 m_{i_1} u_2 m_{i_2} \dots u_k m_{i_k} \bar \gamma$,
with $i_j \! \in \! \{ 0;1 \}$, $u_j \! \in \! U_{i_j}$
and $\bar \gamma \! \in \! B$,}

\pec

or also:
$\gamma = u_1 (m_{i_1} u_2 m_{i_1}^{-1}) \dots
(m_{i_1}m_{i_2} \dots m_{i_{k-1}} u_k m_{i_{k-1}}^{-1} \cdots 
m_{i_2}^{-1} m_{i_1}^{-1})
(m_{i_1}m_{i_2} \dots m_{i_k}) \bar \gamma$,
which can be geometrically interpreted in terms of galleries in a building.

\pec

{\bf 3.E.2 } We can now prove the existence of groups with twin root 
data different from those provided
by Kac-Moody theory.

\pec

T{\sevenrm HEOREM}.--- \it
The group $\Lambda$ satisfies the axioms of a twin root datum for the 
family of root groups
$\{ U_a = u_a({\bf K}_{\iota(a)}) \}_{a \in \Phi}$ above.
As a consequence, there exists a semi-homogeneous Moufang twin tree of
valencies $1 + \mid \! {\bf K}_0 \! \mid$ and $1 + \mid \! {\bf K}_1 
\! \mid$ for any choice
of two fields ${\bf K}_0$ and ${\bf K}_1$.
\rm\pec

{\it Proof}.
The group $\Lambda$ being defined as an amalgam, it acts on a 
semihomogeneous tree $\Delta$ of valencies
$[P_0: B]$ and $[P_1: B]$.
This follows from Bass-Serre theory [S, I.4 Theorem 7].
In view of the decompositions $(\star)$ in 3.E.1, these valencies are
$1 + \mid \! {\bf K}_0 \! \mid$ and $1 + \mid \! {\bf K}_1 \! \mid$.
A fundamental domain for this action is given by the closure of an edge.
The stabilizer of an edge (of a vertex of type 0, resp. of type 1), 
is isomorphic to $B$
(to $P_0$, resp. $P_1$).
In particular, the action is not discrete.
We have the identification of $\Lambda$-sets:

\pec

\centerline{$\Delta \simeq \Lambda/B \sqcup \Lambda/P_0 \sqcup \Lambda/P_1$.}

\pec

The edges (the vertices of type 0, resp. of type 1) are in bijection 
with $\Lambda/B$
(with $\Lambda/P_0$, resp. $\Lambda/P_1$).
We recover the simplicial structure of the tree $\Delta$ thanks to 
the (reversed) inclusion relation on
stabilizers.
Let us set $t:=m_1m_0$.
The set of edges
$\{ t^nB \}_{n \in \hbox{\sevenbf Z}} \sqcup
\{ t^n m_1 B \}_{n \in \hbox{\sevenbf Z}}$
defines a geodesic in the tree: this is our standard apartment ${\Bbb 
A}$, containing the standard edge
$B$.

\pec

We now check the (TRD) axioms for $\Lambda$, as given in 1.A.1, with 
$T$ playing the role of $H$.
We also need root groups.
Recall that the positive roots are described in 3.A.
The root groups indexed by positive roots are the groups of the form
$t^{-n} U_0 t^n$, $t^nm_1 U_0 m_1^{-1}t^{-n}$, $t^n U_1 t^{-n}$ or 
$t^{-n}m_0 U_1 m_0^{-1}t^n$
for $n \geq 0$.
The root groups indexed by negative roots are defined similarly.
In the context of trees, as mentioned in 3.A, a pair of roots $\{ a;b 
\}$ is prenilpotent if $a
\supset b$ or $b \supset a$. Using conjugation by a suitable power of 
$t$ enables us to see prenilpotent
pairs of roots as pairs of positive roots.
Since the group $U_+$ is defined in such a way that the commutation 
relations are trivial for
all prenilpotent pairs, we get axiom (TRD1).
By definition, the group $\Lambda$ is generated by $P_0$ and $P_1$, 
hence by the positive root
groups, by $T$ and by the two root groups indexed by the opposite of 
the simple roots: this is axiom
(TRD4).
Axiom (TRD2) follows from relations in ${\rm SL}_2$ defining the 
elements $m_i(\lambda)$, and from
the definition of the root groups.
 From the definition of $B$, we have that $T$ normalizes the root 
groups, which is the
second half of axiom (TRD0).
By definition of $\Delta$, the simple group $U_i$ is simply 
transitive on the edges whose closure contains
the vertex $P_i$ and different from the standard edge: this proves 
the first assertion of axiom
(TRD0) and the second one in axiom (TRD3).
We are reduced to prove the first assertion of axiom (TRD3), but as 
noted by P. Abramenko
[A \S1 remark 2], this follows from what has already been proved above.
\qed\pec

R{\sevenrm EMARK}. --- A basic consequence of the axioms of a twin 
root datum, is the existence of $BN$-pairs in the group.
For the positive $BN$-pair in $\Lambda$, this can be proved 
concretely as follows.
Since ${\rm Stab}_\Lambda({\Bbb A})$ contains the subgroup
$\langle \ m_0(\lambda_i), m_1(\lambda_1): \lambda_i \! \in \! {\bf 
K}_i^\times \ \rangle$, it is
transitive on the set of edges of ${\Bbb A}$.
Let $\{ e';e'' \}$ be a pair of edges in $\Delta$.
By definition, we have $e'=\gamma' B$ and $e''=\gamma'' B$ for
$\gamma', \gamma'' \! \in \! \Lambda$.
As in the end of the previous subsection, the element 
$\gamma'^{-1}\gamma''$ can be written
$\gamma'^{-1}\gamma'' = (v_1 v_2 \dots v_k) (m_{i_1}m_{i_2} \dots 
m_{i_k}) \bar \gamma$,
with $v_j=m_{i_1}m_{i_2} \dots m_{i_{j-1}} u_j m_{i_{j-1}}^{-1} 
\cdots m_{i_2}^{-1}
m_{i_1}^{-1} \! \in \! U_+$.
Hence we have
$(v_1 v_2 \dots v_k) \gamma'^{-1} \cdot \{ e';e'' \}
=\{ B; m_{i_1}m_{i_2} \dots m_{i_k} B \}$.
In other words, we proved in this paragraph the strong transitivity 
of $\Lambda$ on the tree $\Delta$.
Consequently, the group $\Lambda$ admits a $BN$-pair
$\bigl( B, {\rm Stab}_\Lambda(L) \bigr)$ -- see [Ro1, Theorem 5.2].

\pec

This result (about trees) is stated at a higher level of generality 
and abstraction in [T4].
The down-to-earth computations we used here will be useful in the 
next case, of dimension
two and arbitrarily large rank.

\vfill\eject

\centerline{\bf 4. Twinning Fuchsian buildings}

\pec

The computations of the previous section will now be used to generate 
two-dimensional Moufang twin
buildings of hyperbolic type and arbitrarily large rank.
We will apply the procedure of 2.A, but first we need to examine the 
root system.

\pec

{\bf 4.A } {\it Geometric description of the root system}. --- The 
starting point is the tiling of the hyperbolic plane ${\Bbb H}^2$ by 
regular
right-angled $r$-gons.
We choose such an $r$-gon $R$ and label its edges (open segments)
$\{ E_i \}_{i \in \hbox{\sevenbf Z}/r}$ in a natural cyclic order.
Each $E_i$ supports a geodesic with associated orthogonal reflection $r_i$.
By the Poincar\'e theorem [Mas, IV.H.11], the reflections $\{ r_i 
\}_{i \in \hbox{\sevenbf
Z}/r}$ give rise to a Coxeter system $(W, \{ r_i \}_{i \in 
\hbox{\sevenbf Z}/r})$ and
$W.\overline R = {\Bbb H}^2$ is our tiling. This tiling is a metric 
realization of the
Coxeter complex associated to
$(W, \{ r_i \}_{i \in \hbox{\sevenbf Z}/r})$, which represents the 
spherical facets.
We will use freely the terminology of Coxeter complexes -- see [Ro1 
\S2] and 1.A.1 here.

\pec

L{\sevenrm EMMA}/D{\sevenrm EFINITION}.--- \it
Given a root $a$, the panels on its boundary wall all have the same 
type, and we shall
refer to this as the {\rm type } of the root $a$, and denote it by 
$\iota(a) \! \in \! {\bf Z}/r$.
If $a_i$ is the simple root of type $i$, then we have

\pec

\centerline{${\rm Stab}_W (\partial a_i) =
\langle r_i \rangle \times \langle r_{i-1}, r_{i+1} \rangle 
\simeq{\bf Z}/2 \times D_\infty$
\quad and \quad ${\rm Stab}_W (a_i) = \langle r_{i-1}, r_{i+1} \rangle
\simeq D_\infty$.}

\pec

Hence, we have a bijection between the set of roots of type $i$ and 
$W/\langle r_{i-1}, r_{i+1} \rangle$.
\rm\pec

{\it Proof}.
Let us fix a type $i$ and remark that since we are working with a 
right-angled tiling, the
simple reflections $r_{i-1}$ and $r_{i+1}$ stabilize the wall 
$\partial a_i$, and in fact
the root $a_i$. Moreover $\langle r_{i-1}, r_{i+1} \rangle$ is 
transitive on panels
contained in $\partial a_i$. In particular, $\partial a_i$ is a union 
of closures of panels
of type $i$. This proves our assertion on the type of a root.
Assume now we are given $w \! \in \! {\rm Stab}_W (\partial a_i)$.
The standard panel $E_i \subset \overline R$ of type $i$ is thus sent 
by $w$ on a panel in $\partial a_i$,
which writes $w'E_i$, with $w' \! \in \! \langle r_{i-1}, r_{i+1} 
\rangle$: $w'^{-1}w$ fixes $E_i$, hence
is trivial or equal to $r_i$.
If $w$ stabilizes the root $a_i$, so does $w'$ and we have $w=w'$.
\qed\pec

For a root $a$, we define $\epsilon_a \! \in \! \{ \pm 1 \}$ to be 
$(-1)^m$, where $m$ is the number of
walls of type $\iota(a)$ meeting the interior of the geodesic segment 
from the barycenter of $R$ to the wall
$\partial a$.

\pec

{\bf 4.B } {\it The Borel subgroup. Unipotent radicals of parabolic 
subgroups}. --- The definitions are completely analogous to that of 
the tree case.
For each $i \! \in \! {\bf Z}/r$, we pick a field ${\bf K}_i$.
To each positive root $a$ of type $i$ is attached a copy of the additive
group $({\bf K}_i,+)$.
We denote it by $U_a:= \{ u_a(k): k \! \in \! {\bf K}_i \}$.
As in 1.A.2 we define, for each $w \! \in \! W$, a group $U_w$.
It can be expressed as a product $\prod_{a \in \Phi_{w^{-1}}} U_a$, 
and using the Bruhat ordering of $W$,
these groups $\{ U_w \}_{w \in W}$ form an inductive system with 
$U_w<U_{w'}$ when $w<w'$.

\pec

D{\sevenrm EFINITION}.--- \it
{\rm (i) } The {\rm standard torus } $T$ is the direct product $T:= 
\prod_{i \in {\bf Z}/r} {\bf K}_i^\times$ of
the multiplicative groups of the chosen fields.

{\rm (ii) } The group $U_+$ is the limit of the inductive system 
described above:
$U_+:= \hbox{\rm $\limind$}_{w \in W} U_w$.
\rm\pec

Using the notation of Sect. 2, we view $T$ as the maximal torus
$\Bigl\{ \prod_{i \in {\bf Z}/r} t_i(\lambda_i):
\lambda_i \! \in \! {\bf K}_i^\times \Bigr\}$
of $\prod_{i \in {\bf Z}/r} {\rm SL}_2({\bf K}_i)$.
As in 3.B.1, we make $T$ act on $U_+$ by:

\pec

\hskip 6mm
$t_j(\lambda)u_a(k)t_j(\lambda)^{-1}:=u_a(\lambda^{2\epsilon_a\delta_{j\iota(a)}}k)$.

\pec

In particular, $t_j(\lambda)$ centralizes each root group of type $\neq j$.

\pec

D{\sevenrm EFINITION}.--- \it
The {\rm (positive) Borel subgroup $B$} is the semi-direct product
$B:= T \ltimes U_+$.
\rm\pec

Let us turn now to the construction of \og unipotent radicals\fg.
Call $x_{i,i+1}$ the vertex of type ${i,i+1}$ in the closure 
$\overline R$ of the standard
chamber.
As in 3.A we must classify positive roots according to whether they 
are prenilpotent with
the simple root $a_i$ and/or the simple root $a_{i+1}$.
The point $x_{i,i+1}$ is the vertex of the {\it sector } ${\cal Q}:= 
a_i \cap a_{i+1}$.
Among the four connected components of ${\Bbb H}^2 \setminus 
(\partial a_i \cup \partial a_{i+1})$, this sector is
the one containing $R$.
We have (see the picture below):

\pec

\centerline{
$\overline a_i = \overline {\cal Q} \cup s_{i+1} \overline {\cal Q}$ 
\qquad and \qquad
$\partial a_i = (\overline {\cal Q} \cup s_{i+1} \overline {\cal Q}) \cap
s_i(\overline {\cal Q} \cup s_{i+1} \overline {\cal Q})$. }

\pec

P{\sevenrm ICTURE}. ---

\pec

\centerline{\epsfysize=60mm $$\epsfbox{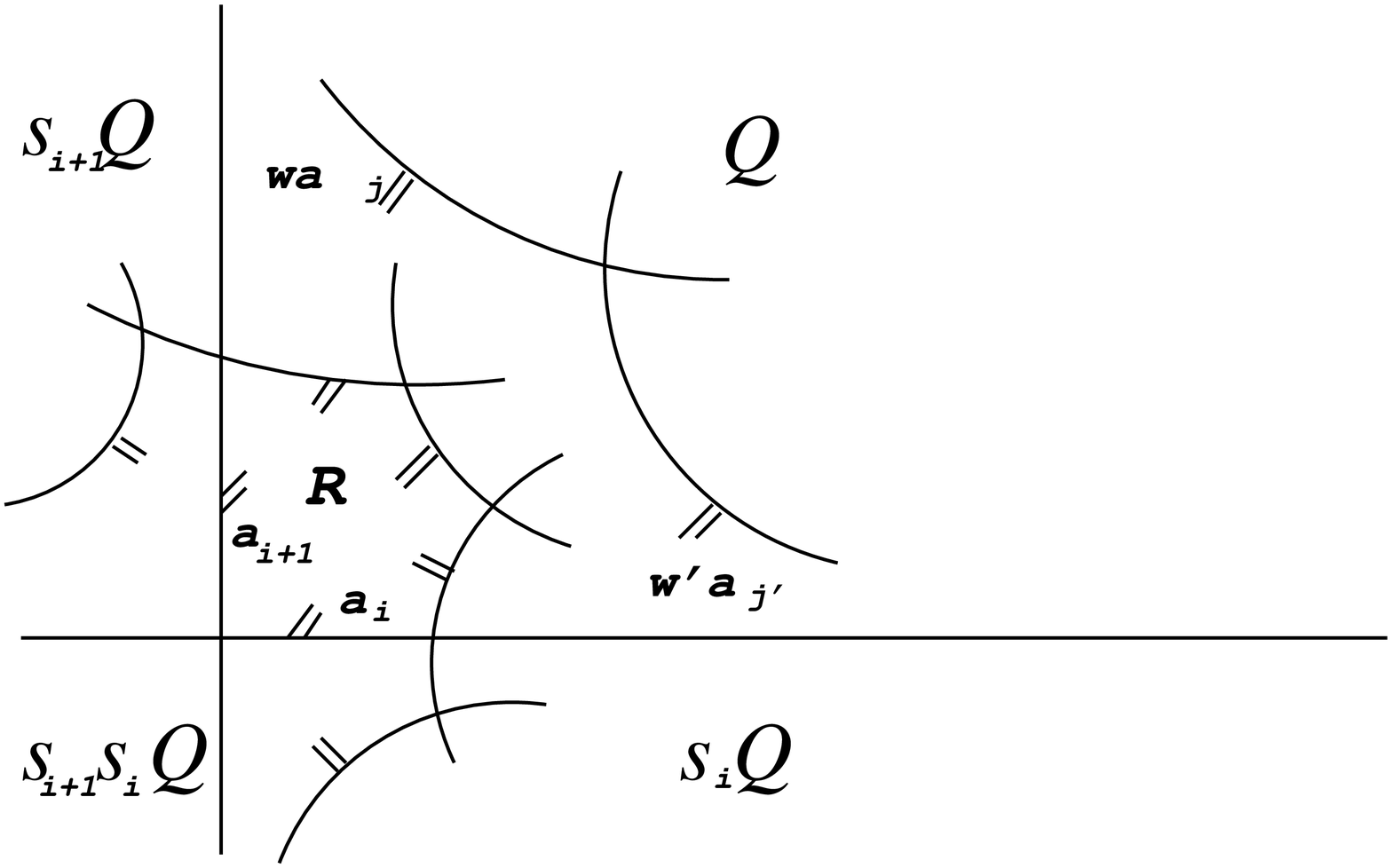}$$}

\pec

D{\sevenrm EFINITION}.--- \it
{\rm (i) } The group $V^{i,i+1}$ is the subgroup of $U_+$ generated 
by all the positive root
groups except $U_i$ and $U_{i+1}$.

{\rm (ii) } The group $U^{i,i+1}$ is the normal closure of 
$V^{i,i+1}$ in $U_+$.
\rm\pec

According to the defining relations of $U_+$, the group $U^{i,i+1}$ 
admits the following infinite set of
generators (satisfying the obvious induced relations), where $a$ is a 
positive root.

\pec

\centerline{$u_i(v)u_{i+1}(t)u_a(k)u_{i+1}(t)^{-1}u_i(v)^{-1}$ for 
$\partial a \subset {\cal Q}$.}

\pec

These generators will be referred to as the generators {\it of the first type}.

\pec

\centerline{$u_i(v)u_a(k)u_i(v)^{-1}$ \quad $\left\{\matrix{
\hbox{\rm for $\partial a \subset s_{i+1} {\cal Q}$,} \hfill \cr
\hbox{\rm or if $\partial a \cap \partial a_{i+1}$ is a point in 
$s_{i+1} \overline{\cal Q}$.}
}\right.$}

\pec

These generators will be referred to as the generators {\it of the 
second type}.

\pec

\centerline{$u_{i+1}(t)u_a(k)u_{i+1}(t)^{-1}$ \quad $\left\{\matrix{
\hbox{\rm for $\partial a \subset s_i {\cal Q}$,} \hfill \cr
\hbox{\rm or if $\partial a \cap \partial a_i$ is a point in $s_i 
\overline{\cal Q}$.}
}\right.$}

\pec

These generators will be referred to as the generators {\it of the third type}.

\pec

\centerline{$u_a(k)$ \quad $\left\{\matrix{
\hbox{\rm for $\partial a \subset s_i s_{i+1} {\cal Q}$,} \hfill \cr
\hbox{\rm if $\partial a \cap \partial a_i$ is a point in $s_i 
s_{i+1} \overline{\cal Q}$,} \hfill \cr
\hbox{\rm or if $\partial a \cap \partial a_{i+1}$ is a point in $s_i 
s_{i+1} \overline{\cal Q}$.}
}\right.$}

\pec

These generators will be referred to as the generators {\it of the 
fourth type}.

\pec

R{\sevenrm EMARK}. --- The parameters $v$, $t$ and $k$ are in the 
fields ${\bf K}_i$, ${\bf K}_{i+1}$ and
${\bf K}_{\iota(a)}$, respectively.

\pec

{\bf 4.C } {\it Actions of Levi factors on unipotent radicals}. --- 
We first define the Levi factors.

\pec

D{\sevenrm EFINITION}.--- \it
The {\rm (standard) Levi factor of type $i,i+1$ } is the direct product

\pec

\centerline{$L_{i,i+1}:=
{\rm SL}_2({\bf K}_i) \times {\rm SL}_2({\bf K}_{i+1}) \times 
\prod_{j \neq i,i+1}
t_j({\bf K}_j^\times)$.}
\rm\pec

By definition, this group contains the torus $T$; it is generated by 
the toric elements $t_j(\lambda_j)$
($j \neq i, i+1$), the unipotent elements $u_i(r)$ ($r \! \in \! {\bf 
K}_i$) and
$u_{i+1}(s)$ ($s \! \in \! {\bf K}_{i+1}$), and by the elements 
$m_i(\lambda_i)$
($\lambda_i \! \in \! {\bf K}_i^\times$) and $m_{i+1}(\lambda_{i+1})$
($\lambda_{i+1} \! \in \! {\bf K}_{i+1}^\times$).

\pec

Let us turn now to the action of such a group $L_{i,i+1}$ on $U^{i,i+1}$.
We shall specify that:

\pec

-- an element $u_i(r)$ centralizes the generators of the third and 
fourth types, and changes the
conjugating element $u_i(v)$ into $u_i(r+v)$ for the generators of 
the first and second types, as
in (3C1) and (3C2).

\pec

-- an element $u_{i+1}(s)$ centralizes the generators of the second 
and fourth types, and
changes the conjugating element $u_{i+1}(t)$ into $u_{i+1}(s+t)$ for 
the generators of the first
and third types, as in (3C1) and (3C2).

\pec

-- an element $t_j(\lambda_j)$ centralizes all the root groups 
indexed by roots of type $\neq j$.
On a root group $U_a$ with $\iota(a)=j$, it acts by multiplication of 
the additive parameter by
$\lambda_j^{2\epsilon_a}$, as in (3B1).

\pec

Using (3C3) and (3C4) for trees, the action of the elements 
$m_i(\lambda_i)$ and
$m_{i+1}(\lambda_{i+1})$ is defined as follows.
For generators of the first type and when $v \neq 0$, we set:

\pec

\hskip 6mm
$m_i(\lambda) \bigl( u_i(v)u_{i+1}(t)u_a(k)u_{i+1}(t)^{-1}u_i(v)^{-1} \bigr)
m_i(\lambda)^{-1}$

\hskip 12mm
$\displaystyle
:= u_i({-\lambda^{2} \over v}) u_{i+1}(t)
u_a \bigl( ({-\lambda \over v})^{2\epsilon_a\delta_{i\iota(a)}}k 
\bigr) u_{i+1}(t)^{-1}
u_i({-\lambda^{2} \over v})^{-1}$;

\pec

whereas when $v=0$, we set:

\pec

\hskip 6mm
$m_i(\lambda) \bigl( u_{i+1}(t)u_a(k)u_{i+1}(t)^{-1} \bigr) m_i(\lambda)^{-1}:=
u_{i+1}(t) u_{s_ia} (\lambda^{-2\epsilon_a\delta_{i\iota(a)}}k) 
u_{i+1}(t)^{-1}$.
\pec

For generators of the second type and when $v \neq 0$, we set:

\pec

\hskip 6mm
$\displaystyle
m_i(\lambda) \bigl( u_i(v)u_a(k)u_i(v)^{-1} \bigr) m_i(\lambda)^{-1}
:= u_i({-\lambda^{2} \over v}) u_a \bigl( ({-\lambda \over 
v})^{2\epsilon_a\delta_{i\iota(a)}}k \bigr)
u_i({-\lambda^{2} \over v})^{-1}$;

\pec

whereas when $v=0$, we set:

\pec

\hskip 6mm
$m_i(\lambda)u_a(k)m_i(\lambda)^{-1}:= u_{s_ia}
(\lambda^{-2\epsilon_a\delta_{i\iota(a)}}k)$.

\pec

For generators of the third type, we set:

\pec

\hskip 6mm
$m_i(\lambda) \bigl( u_{i+1}(t)u_a(k)u_{i+1}(t)^{-1} \bigr) m_i(\lambda)^{-1}:=
u_{i+1}(t)u_{s_ia}(\lambda^{-2\epsilon_a\delta_{i\iota(a)}}k)u_{i+1}(t)$.

\pec

For generators of the fourth type, we set:

\pec

\hskip 6mm $m_i(\lambda)u_a(k)m_i(\lambda)^{-1}:=
u_{s_ia}(\lambda^{-2\epsilon_a\delta_{i\iota(a)}}k)$.

\pec

R{\sevenrm EMARK}. ---
As in 3.C, $\delta_{i\iota(a)}$ in the exponents 
${}^{2\epsilon_a\delta_{i\iota(a)}}$
means that the element $m_i(\lambda)$ induces a multiplication of the 
additive parameter $k$ in
$u_a(k)$ by a factor $\lambda^{-2\epsilon_a}$ or
$\displaystyle ({-\lambda \over v})^{2\epsilon_a}$ only if the type 
$\iota(a)$ of the root $a$ is $i$;
otherwise, it doesn't change $k$.

\pec

This defines the action of the factor ${\rm SL}_2({\bf K}_i)$; that 
of ${\rm SL}_2({\bf K}_{i+1})$ is
defined {\it mutatis mutandis}.

\pec

{\bf 4.D } {\it Checking the product relation}. ---
According to the definition of the Levi factor of type $i,i+1$ as a 
direct product, it is enough to check
separately the defining relations for each toric factor $t_j({\bf 
K}_j^\times)$ ($j \neq i,i+1$) and each of
both factors ${\rm SL}_2({\bf K}_i)$ and ${\rm SL}_2({\bf K}_{i+1})$, 
whose actions obviously commute
to one another.

\pec

Let us consider the ${\rm SL}_2({\bf K}_i)$-action, and check the 
product relation
$g(hvh^{-1})g^{-1}=(gh)v(gh)^{-1}$.
We consider a generator $v=u_i(v)u_a(k)u_i(v)^{-1}$ of $U^{i,i+1}$ of 
the second type.
Then the root $a$ is not prenilpotent with $a_i$ and, in view of the 
definition of the actions in 4.C, the
computations to check the product relation reduce to the computation 
made in 3.D for a generator of the
first type in the sense of trees (3.D).
Similarly, the computations to check the product relation for a 
generator of the fourth type in the
sense of Fuchsian buildings (4.B) reduce to the computations made in 
3.D for a generator of the second
type in the sense of trees (3.D).
Finally, by the defintions in 4.C, the actions of $u_i(r)$ and of 
$m_i(\lambda_i)$ commute with
conjugation by the elements $u_{i+1}(t)$.
So, up to the conjugation by $u_{i+1}(t)$, checking the product 
relation for a generator of the
third (resp. first) type in the sense of Fuchsian buildings amounts 
to checking the product relation for
a generator of the second (resp. first) type in the sense of trees.

\pec

Consequently, we can introduce the following definition.

\pec

D{\sevenrm EFINITION}.--- \it
{\rm (i) } The {\rm (standard) parabolic subgroup of type $i,i+1$ } 
is the semidirect product

\pec

\centerline{$P_{i,i+1}:= L_{i,i+1} \ltimes U^{i,i+1}$.}

\pec

{\rm (ii) } The {\rm (standard) parabolic subgroup of type $i$ } is 
the subgroup $P_i$ of
$P_{i,i+1}$ generated by $L_i:= \langle T,{\rm SL}_2({\bf K}_i) 
\rangle$ and $U_+$.
\rm\pec

The groups $P_i$ have a Levi decomposition
$P_i = L_i \ltimes U^i$, where $U^i$ is the normal closure in $U_+$ of all the
positive root groups except $U_i$, and $L_i$ is of course $\langle 
T,{\rm SL}_2({\bf K}_i) \rangle$.
Besides for each $i \! \in \! {\bf Z}/r$, we have chains of inclusions:
$B \hookrightarrow P_i \hookrightarrow P_{i,i+1}$
and $B \hookrightarrow P_i \hookrightarrow P_{i,i-1}$.
This gives rise to an inductive system of group homomorphisms for 
which we will have a geometric
interpretation in the next subsection.
As for trees, we will amalgamate the parabolic subgroups in order to 
define a group $\Lambda$ endowed with a twin root datum.

\pec

{\bf 4.E } {\it Group combinatorics}. --- The verification of the 
axioms of a twin root datum for $\Lambda$ will of course follow the 
lines of 3.E.
Nevertheless, we first need to present some geometric notions which 
are used to prove the existence of a
building acted upon by $\Lambda$.

\pec

{\bf 4.E.1 } Let us recall some facts about complexes of groups [BH 
III.${\cal C}$], as well
as an application to hyperbolic buildings due to D. Gaboriau and F. 
Paulin [GP].
The use of this theory should not be surprising, because it is a
higher-dimensional generalization of Bass-Serre theory about group 
actions on trees, and we used
Bass-Serre theory in the proof of 3.E.2.
We will only use the simpler notion of {\it polytope of groups}, 
namely the datum of a compact convex
polytope $R$, of a group $G_\sigma$ for each face $\sigma$ of $R$ and 
of an injective group
homomorphism $G_\sigma \hookrightarrow G_\tau$ for each inclusion of 
faces $\tau \subset \bar
\sigma$.
We require the commutativity of the diagram of group homomorphisms 
given by all the
(reversed) inclusions of faces.
This diagram is an inductive system indexed by the barycentric 
subdivision of $R$.
The inductive limit is the {\it fundamental group } of the polytope 
of groups [BH III.${\cal C}$.3].

\pec

The connection with 4.D is that $(R,\{B\hookrightarrow P_i\hookrightarrow
P_{i,i+1}\}_{i\in{\bf Z}/r})$ is a polytope of groups.
The group $B$ (resp. $P_i$, resp. $P_{i,i+1}$) is attached to the 
2-dimensional cell of $R$ (resp.
to the edge of type $\{ i \}$, resp. to the vertex of type $\{i,i+1\}$).

\pec

D{\sevenrm EFINITION}.--- \it
The group $\Lambda$ is the limit of the finite inductive system 
described in {\rm 4.D}.
In other words, it is the fundamental group of the above complex of groups.
\rm\pec

In the case of graphs, Bass-Serre theory provides a tree on which the 
fundamental group acts, with the
groups $G_\sigma$ as prescribed stabilizers.
In the case of polytopes of groups, this kind of existence result is 
not immediate, because non-positive curvature
arguments come into play [BH, 4.17].
Under non-positive curvature assumptions, we know that $\Lambda$ acts 
on a simply connected cell complex
with the groups $G_\sigma$ as prescribed stabilizers.
In our situation, this will be even better, since we can apply the 
result of D. Gaboriau and F. Paulin
alluded to above.
We follow their terminology:

\pec

D{\sevenrm EFINITION}.--- \it
If $R$ is a hyperbolic polyhedron providing a Poincar\'e tiling of 
the hyperbolic space
${\Bbb H}^n$, a {\rm hyperbolic building of type $R$ } is a piecewise 
polyhedral cell
complex, covered by a family of subcomplexes -- the {\rm apartments} 
-- all isomorphic to this
tiling and satisfying the following incidence axioms.

\pec

{\rm (i)} Two points are always contained in an apartment.

{\rm (ii)} Two apartments are isomorphic by a polyhedral arrow fixing their
intersection.

\pec

A building whose apartments are tilings of the hyperbolic plane 
${\Bbb H}^2$ will be
called {\rm Fuchsian}.
\rm\pec

R{\sevenrm EMARK}. --- A hyperbolic building carries a natural 
CAT$(-1)$ metric [GP, Proposition 1.5].

\pec

Let us go back to our construction.
According to [GP, Theorem 0.1], there exists a hyperbolic building 
$\Delta$ of type $R$ such that we
have the identification of $\Lambda$-sets:

\pec

\centerline{$\displaystyle \Delta \simeq \Lambda/B \sqcup
\bigsqcup_{i \in {\bf Z}/r} \Lambda/P_i \sqcup
\bigsqcup_{i \in {\bf Z}/r} \Lambda/P_{i,i+1}$.}

\pec

The building structure is given by the (reversed) inclusion relation 
on stabilizers.
For $\gamma \! \in \! \Lambda$, a translate $\gamma B$ (resp.
$\gamma P_i$, $\gamma P_{i,i+1}$) is a {\it chamber } (resp. an {\it 
edge of type $i$},
a {\it vertex of type $i,i+1$}).
A fundamental domain for this action is given by the closure of $R$.
Moreover it is shown in [Bou3] that given $R$ and
$\underline q:= \{ q_i \}_{1 \leq i \leq r}$ a sequence of integers 
$\geq 2$, there exists
a unique Fuchsian building $I_{r,1+\underline q}$ with apartments 
isomorphic to the tiling
of ${\Bbb H}^2$ by $R$, and such that the link at any vertex of type 
$i,i+1$ is the complete
bipartite graph of parameters $(1+q_i,1+q_{i+1})$. Recall that the 
{\it link } at a point is a
sufficiently small sphere centered at this point; in our 
two-dimensional context, it is seen
as a graph. Uniqueness implies that when we choose finite fields for 
our construction, the
building $\Delta$ above is
$I_{r, \underline q+1}$, for $\underline q$ a sequence of prime powers.
When the $q_i$'s are all equal to a given prime power $q$, the 
building comes from a (non-unique) Kac-Moody group;
we denote it by $I_{r,1+q}$.

\pec

R{\sevenrm EMARK}. --- Let $\Lambda$ be a Kac-Moody group over ${\bf 
F}_q$ whose Weyl group is the group $W$ associated to the
hyperbolic tiling of ${\Bbb H}^2$ we are considering.
Such a group exists: choose any generalized Cartan matrix 
$A=[A_{i,j}]_{i,j \in {\bf Z}/r}$
such that $A_{i,i}=2$, $A_{i,i+1}=0$ and $A_{i,j} A_{j,i} \geq 4$ for 
$j \neq i, i+1$.
In view of its combinatorial structure (positive $BN$-pair) and of 
Tits' amalgam theorem [T1 \S14],
$\Lambda$ is the limit of the inductive system given by the 
inclusions of spherical parabolic subgroups
of rank $\leq 2$.
That is, the fundamental group of a complex of groups with the 
barycentric subdivision of a regular
right-angled $r$-gon as indexing polytope.
In the construction we give, in which the base field can vary from 
one panel of the base chamber to another, the parabolic
subgroups are designed to have a similar Levi decomposition to these 
Kac-Moody examples.

\pec

{\bf 4.E.2 } Here is now our main constructive result about twinnings 
in the two-dimensional
case of Fuchsian buildings.

\pec

T{\sevenrm HEOREM}.--- \it
The group $\Lambda$ defined in {\rm 4.E.1 } satisfies the axioms of a 
twin root datum for the above family of root groups
$\{ U_a = u_a({\bf K}_{\iota(a)}) \}_{a \in \Phi}$.
In particular, given any regular right-angled $r$-gon $R$, and a set 
$\underline q:= \{ q_i \}_{1 \leq i \leq r}$ of prime powers,
the above defined Fuchsian building $I_{r,1+\underline q}$ belongs to a Moufang
twinning.
\rm\pec

{\it Proof}.
The aim of the first paragraph of the proof of 3.E.2 was to make a bit more 
explicit the use of Bass-Serre
theory.
In the higher-dimensional case of Fuchsian buildings, the analogous 
work was done in 4.E.1, where we
presented what we need from H\ae fliger's theory of complexes of groups.
The verification of the axioms of a twin root datum is then the same 
as in 3.E.2.
The final statement follows from the uniqueness of the Fuchsian 
building [Bou3] with the given local data.
\qed\pec

R{\sevenrm EMARK}. --- As in remark 3.E.2, we can see directly that
$\bigl( B, {\rm Stab}_\Lambda({\Bbb A}) \bigr)$ is a $BN$-pair in $\Lambda$.

\vfill\eject

\centerline{\bf 5. Non-linearities}

\pec

A Kac-Moody group defined over a finite ground-field ${\bf F}_q$ has 
no reason to be linear over
${\bf F}_q$ merely because it is defined by a presentation involving 
subgroups that are linear over
${\bf F}_q$.
Still, in case it {\it is} linear over a field, then the 
characteristic of the latter field must be that of its ground field
(5.A).
In Sect. 3 and 4 we constructed groups of Kac-Moody type with more 
than one ground field, and in 5.B
we show that if they involve fields with two different 
characteristics, then they cannot be linear.

\pec

{\bf 5.A } {\it Negative results for Kac-Moody groups}. ---
In this subsection, we use the notation of 1.C, because we are 
specifically working with a topological Kac-Moody group $G$.
We denote by $(G,N,U_c,\Gamma,T,S)$ its refined Tits system (1.C.2).
The appropriate axioms, given in 1.A.5, can be proved for the affine 
$BN$-pairs associated to Chevalley
groups over local fields of equal characteristic to their residue 
fields, but not in the case of local fields of
characteristic 0.  The following result can be found in [R\'e3].

\pec

P{\sevenrm ROPOSITION}.--- \it
Assume $\Lambda$ is an infinite Kac-Moody group over ${\bf F}_q$.
Then the groups $\Lambda$ and $G$ cannot be linear over any field of 
characteristic different
from $p$.
\qed\rm\pec

{\bf 5.B } {\it The wider context of twin root data}. --- This is the 
point where Moufang twinnings not
arising from Kac-Moody theory provide lattices of hyperbolic 
buildings which are not linear at all.
The result below should be seen as a consequence of the constructive 
theorem 4.E.2.

\pec

T{\sevenrm HEOREM}.--- \it
Let $r\geq 5$ be an integer and $\{ {\bf K}_i \}_{i \in {\bf Z}/r}$ 
be a family of $r$ fields.
Assume there exist two distinct indices in ${{\bf Z}/r}$ such that 
the corresponding fields have different
positive characteristics.
Let $\Lambda$ be the group defined as in {\rm 4.E.2} from these fields, and
$\Gamma$ be the fixator of a chamber in the associated Moufang 
Fuchsian twinning.
Then, any group homomorphism

\pec

\centerline{$\displaystyle \rho: \Gamma \to \prod_{\alpha \in A} 
{\Bbb G}_\alpha({\bf F}_\alpha)$}

\pec

has infinite kernel, whenever the index set $A$ is finite and ${\Bbb 
G}_\alpha$ is a linear algebraic group
defined over the field ${\bf F}_\alpha$ for each $\alpha \! \in \! A$.
\rm\pec

R{\sevenrm EMARK}. ---
In particular, there exist twin root data for Fuchsian twinnings that 
admit non-linear lattices
(since for large enough thicknesses, a chamber-fixator is a lattice 
in the full automorphism group of the
building of sign opposite the chamber).

\pec

{\it Proof}.
Since any two chamber-fixators are isomorphic, we can deal with the 
standard positive chamber $R$ of 4.A.
Then, the corresponding chamber-fixator is the Borel subgroup $B=T 
\ltimes U_+$ of subsect. 4.B, and
for the proof we shall use the notation $B$ rather than $\Gamma$.
As before, ${\Bbb A}$ is the apartment stabilized by the group $N$ 
generated by the elements
$t_i(\lambda)$ and $m_i(\lambda)$ when $i$ ranges over ${\bf Z}/r$ 
and $\lambda$ ranges over
${\bf K}_i^\times$.

\pec

The element $m_i(1)\circ m_{i+2}(1)$ is a hyperbolic translation 
$t_{i+1}$ in ${\Bbb A}$ along the wall
$\partial a_{i+1}$.
The sequence $\{ a_i(n):=t_{i+1}^n.a_i \}_{n \geq 0}$ of positive 
roots has the property that
$a_i(n+1) \supset a_i(n)$.
Moreover for any $n,n' \geq 0$, the root groups $U_{a_i(n)}$ and 
$U_{a_i(n')}$ commute with one another,
so the group $V_i: = \langle U_{a_i(n)} : n \geq 0 \rangle$ is 
isomorphic to the additive group of the polynomial ring
${\bf K}_i[X]$.
In particular, $V_i$ is an infinite group of exponent $p_i:={\rm 
char}({\bf K}_i)$.

\pec

The assumption $r \geq 5$ allows us to choose fields ${\bf K}_i$ and 
${\bf K}_j$ with characteristics
$p_i \neq p_j$ such that the panels of $R$ having types $i$ and $j$ 
support parallel walls.
This defines as above an infinite group $V_j$ of exponent $p_j$.

\pec

P{\sevenrm ICTURE}. ---

\pec

\centerline{\epsfysize=50mm $$\epsfbox{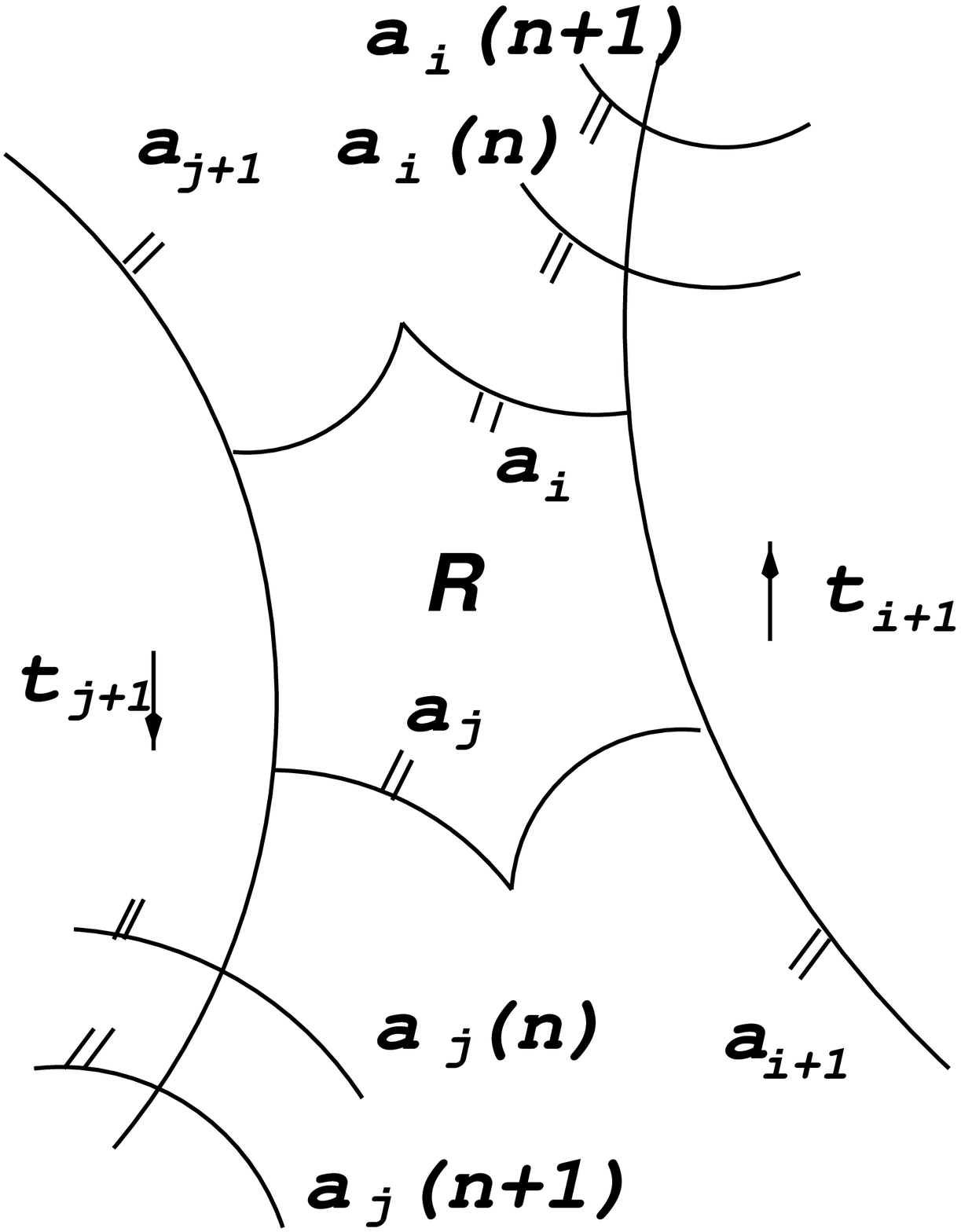}$$}

\pec

We set $A_i: = \{ \alpha \! \in \! A: $ char$({\bf F}_\alpha)=p_i\}$ and
$A_j:= \{ \alpha \! \in \! A: $ char$({\bf F}_\alpha)=p_j\}$.
For each $\alpha \! \in \! A$, we denote by
${\rm pr}_\alpha: \prod_{\alpha \in A} {\Bbb G}_\alpha({\bf 
F}_\alpha) \twoheadrightarrow
{\Bbb G}_\alpha({\bf F}_\alpha)$ the natural projection.
By [Mar, VIII.3.7] for each $\alpha \not\in A_i$ the group
$({\rm pr}_\alpha \circ \rho)(V_i)$ is finite, so there is a finite 
index normal subgroup
$N_i \triangleleft V_i$ such that $\prod_{\alpha \not\in A_i} ({\rm 
pr}_\alpha \circ \rho)(N_i)$ is
trivial, and similarly for $j$ replacing $i$.

\pec

The subgroup generated by the positive and negative root groups of 
type $i$ and $j$ satisfies the axioms of a twin root
datum with infinite dihedral Weyl group, and by applying [KP, 
proposition 3.5 (c)] to the corresponding refined Tits
system (see 1.A.5), one has a free product decomposition of the 
subgroup generated by the positive root
groups of type $i$ and $j$.
This shows that $V_i*V_j$ injects in $B$, and if we pick $v \! \in \! 
N_i \setminus \{1\}$ and $v' \! \in
\! N_j \setminus \{1\}$, then $vv'$ has infinite order.
But since $v \! \in \! N_i$ and $v' \! \in \! N_j$ the images 
$\rho(v)$ and $\rho(v')$ commute with one another and
both have finite order.
We have found a subgroup of $B$ isomorphic to ${\bf Z}$ and with 
finite image under $\rho$.
This proves the theorem.
\qed\pec

R{\sevenrm EMARK}. --- This result shows that closures of arbitrary 
groups with twin root
data are too wide a framework for linearity problems.
Hence, the last linearity problem to be solved is that of remark 5.A, 
the equal characteristic
case for Kac-Moody groups over finite fields.

\pec

{\bf 5.C } {\it Analogy with trees}. ---
In this final subsection we specialize to the case of Kac-Moody 
groups, where all the fields ${\bf K}_i$ are
equal to the finite field of $q$ elements.

\pec

P{\sevenrm ROPOSITION}.--- \it
For any prime power $q \geq 3$, there exists a Kac-Moody group $\Lambda$ over
${\bf F}_q$ whose building is isomorphic to the right-angled Fuchsian 
building $I_{r,1+q}$ and
such that its natural image in ${\rm Aut}(I_{r,1+q})$ contains a 
uniform lattice $\Gamma_{r,1+q}$
abstractly defined by the presentation

\pec

\centerline{$\Gamma_{r,1+q}=
\langle \, \{ \gamma_i \}_{i \in {\bf Z}/r} \,  : \, \gamma_i^{q+1}=1$ and
$[\gamma_i,\gamma_{i+1}]=1 \, \rangle$.}
\rm\pec

R{\sevenrm EMARK}. ---
The above uniform lattices $\Gamma_{r,1+q}$ (as well as the buildings 
$I_{r,1+q}$ themselves) were
introduced by M. Bourdon [Bou2].

\pec

For the proof below, it is good to keep in mind the basic facts on 
Kac-Moody groups recalled in 2.E.1.

\pec

{\it Proof}.
We fix the finite field ${\bf F}_q$ and consider the generalized 
Cartan matrix $A$ indexed by ${\bf Z}/r$
and defined as in 2.E.2 by $A_{i, i}=2$, $A_{i,i+1}=0$ and $A_{i, 
j}=1-q$ for $j \neq i, i\pm1$.
In order to define a full Kac-Moody root datum, we need ${\bf 
Z}$-lattices. In our case, we start by the
lattice of cocharacters, which we define as
$X_*:=\bigoplus_{j \in {\bf Z}/r} {\bf Z}h_j \oplus {\bf Z}\xi_j$.
The ${\bf Z}$-lattice $X^*$ is by definition the ${\bf Z}$-dual of 
$X_*$, and the elements $a_i$ are those
in $X^*$ which satisfy $a_i(h_j)=A_{ij}$ and 
$a_i(\xi_j)=-\delta_{ij}$ for all $i$ and $j$ in ${\bf Z}/r$.
The Kac-Moody root datum
$\bigl({\bf Z}/r, A, X^*, \{a_i\}_{i \in {{\bf Z}/r}}, \{h_i\}_{i \in 
{{\bf Z}/r}} \bigr)$
and the finite field ${\bf F}_q$ define a Kac-Moody group $\Lambda$ 
with maximal split torus
$T={\rm Hom}({\bf F}_q[X^*],{\bf F}_q)$.
For any $j \! \in \! {\bf Z}/r$, let us denote by $\Lambda_j$ the 
subgroup generated by $T$ and the root
groups $U_{\pm a_j}$ indexed by the simple root $a_j$ and its opposite.
This group is the Kac-Moody group with Kac-Moody \og rank-one sub-root datum\fg
$\bigl( j, [2], X^*, a_j, h_j \bigr)$ and is isomorphic to
${\rm GL}_2({\bf F}_q) \times ({\bf F}_q^\times)^{2(r-1)}$.
Since it is the Levi factor of the standard parabolic subgroup of 
type $\{ j \}$ [R\'e1, 6.2.2], we deduce
that the Levi factor of a spherical parabolic subgroup of type
$J$ is isomorphic to $({\bf F}_q^\times)^{2r}$, to the direct product 
${\rm GL}_2({\bf
F}_q) \times ({\bf F}_q^\times)^{2(r-1)}$ or to $\bigl( {\rm 
GL}_2({\bf F}_q) \times {\rm GL}_2({\bf F}_q) \bigr)
\times ({\bf F}_q^\times)^{2(r-2)}$, according to whether $\mid \! J 
\! \mid$ equals 0, 1 or 2.

\pec

By [R\'e1, 9.6.2], the center $Z(\Lambda)$ of the group $\Lambda$ is 
the subgroup of $T$ whose group
of cocharacters is $\bigoplus_{i \in {\bf Z}/r} {\bf Z}(h_i+2\xi_i)$.
It is the subgroup of $T$ which centralizes all root groups, and by 
the argument of Lemma 1.B.1,
$Z(\Lambda)$ is also the kernel of the homomorphism $\varphi: \Lambda \to {\rm
Aut}(I_{r,1+q})$ attached to the action of $\Lambda$ on any of its 
two buildings.
We have $Z(\Lambda) \simeq ({\bf F}_q^\times)^r$, and for each $i \! 
\in \! {\bf Z}/r$, the subgroup
$Z(\Lambda)$ intersects the factor isomorphic to ${\rm GL}_2({\bf F}_q)$ of
$\Lambda_i $ along its torus of scalar matrices.
Consequently, the Levi factors of the parabolic subgroups of 
$\varphi(\Lambda)$ are
isomorphic to
$({\bf F}_q^\times)^r$, to ${\rm PGL}_2({\bf F}_q) \times ({\bf 
F}_q^\times)^{r-1}$
and to $\bigl( {\rm PGL}_2({\bf F}_q) \times {\rm PGL}_2({\bf F}_q) 
\bigr) \times ({\bf F}_q^\times)^{r-2}$.

\pec

On the one hand, as noted in 4.E.2, it is a consequence of Tits' 
theorem [T1 \S 14]
that $\varphi(\Lambda)$ is the fundamental group of a complex of 
groups indexed by (the barycentric subdivision
of) a regular right-angled $r$-gon $R$.
The involved groups are the standard spherical parabolic subgroups of 
rank $\leq 2$, seen as facet
fixators. On the other hand, the group
$\Gamma_{r,1+q}$ is also the fundamental group of a complex of groups indexed
by the same $r$-gon $R$.
The group attached to the edge of type $i$ (resp. the vertex of type 
$i,i+1$) is ${\bf Z}/(q+1)$ (resp.
${\bf Z}/(q+1) \times {\bf Z}/(q+1)$).
The group attached to $R$ itself is $\{ 1 \}$.
Hence, in order to see $\Gamma_{r,1+q}$ as a subgroup of 
$\varphi(\Lambda)$, it is enough to prove the existence of a
morphism from the latter inductive system of finite groups to the 
former inductive system of parabolic
subgroups, since taking the limit will then show that 
$\Gamma_{r,1+q}<\varphi(\Lambda)$.

\pec

But as in [Ch, Proposition 5], we can see ${\bf Z}/(q+1)$ as a subgroup of
${\rm PGL}_2({\bf F}_q)$ which is simply transitive on the chambers 
of the corresponding building ${\Bbb P}^1{\bf F}_q$.
The inclusions
$\{ 1 \}<({\bf F}_q^\times)^r$, ${\bf Z}/(q+1)<{\rm PGL}_2({\bf F}_q) 
\times ({\bf F}_q^\times)^{r-1}$
and ${\bf Z}/(q+1) \times {\bf Z}/(q+1)<\bigl( {\rm PGL}_2({\bf F}_q) 
\times {\rm PGL}_2({\bf F}_q) \bigr) \times ({\bf F}_q^\times)^{r-2}$
then provide the morphism of inductive systems we are looking for.
\qed\pec

R{\sevenrm EMARKS}. ---
1) The group $\Gamma_{r,1+q}$ is a uniform lattice of the building
$I_{r,1+q}$. Indeed, by the very definition of a fundamental group of 
a complex of groups, attaching $\{ 1
\}$ to the 2-cell $R$ implies that $\Gamma_{r,1+q}$ is simply 
transitive on the chambers of the buildings.
Finiteness of all facet fixators implies discreteness.

2) When the number $r$ of edges is even and satisfies
$\sin (\pi/r)<1/\sqrt{q+1}$, the group $\Gamma_{r,1+q}$ can be 
embedded as a convex cocompact
discrete group of isometries of the real hyperbolic space of 
dimension $2q$ [Bou1].

\pec

C{\sevenrm OROLLARY}.--- \it
Whenever $q$ is large enough and $r$ is even and large enough with 
respect to $q$, the topological
Kac-Moody group $G$ defined from the above group $\Lambda$ as in {\rm 
1.B.2}, has virtually pro-$p$ maximal compact
subgroups, and contains both uniform lattices which are linear in 
characteristic $0$ and non-uniform lattices
such that the only possible characteristic of linearity is $p$ 
$($dividing $q)$.
\rm\pec

{\it Proof}. By the non-positive curvature property and the 
Bruhat-Tits fixed point lemma, the maximal compact subgroups are the
vertex-fixators, hence the claim on virtual pro-$p$-ness of these 
groups by 1.C.2 (ii).
Assume $q$ large enough, so that the spherical parabolic subgroups of 
positive (resp. negative) sign
are lattices of the negative (resp. positive) building -- see 1.B.3 and 1.C.1.
It then remains to assume that $r$ (even) satisfies $\sin 
(\pi/r)<1/\sqrt{q+1}$ and apply
Bourdon's embedding result quoted in remark 2) above.
\qed\gec

\centerline{\bf References}

\pec

[A] P. A{\sevenrm BRAMENKO},
Twin buildings and applications to $S$-arithmetic groups.
Lecture Notes in Mathematics {\bf 1641}, Springer, 1997

[BH] M. B{\sevenrm RIDSON}, A. H{\sevenrm \AE FLIGER},
Metric spaces of non-positive curvature.
Grund. der Math. Wiss. {\bf 319}, Springer, 1999

[Bou1] M. B{\sevenrm OURDON},
{\it Sur la dimension de Hausdorff de l'ensemble limite d'une famille 
de sous-groupes convexes
cocompacts}.
C. R. Acad. Sci. Paris {\bf 325} (1997) 1097-1100

[Bou2] M. B{\sevenrm OURDON},
{\it Immeubles hyperboliques, dimension conforme et rigidit\'e de Mostow}.
GAFA {\bf 7} (1997) 245-268

[Bou3] M. B{\sevenrm OURDON},
{\it Sur les immeubles fuchsiens et leur type de quasi-isom\'etrie}.
Erg. Th. and Dynam. Sys. {\bf 20 } (2000), 343-364

[BM] M. B{\sevenrm URGER}, S. M{\sevenrm OZES},
{\it ${\rm CAT}(-1)$-spaces, divergence groups and their commensurators}.
J. of AMS {\bf 9} (1996) 57-93

[BP] M. B{\sevenrm OURDON}, H. P{\sevenrm AJOT},
{\it Rigidity of quasi-isometries for some hyperbolic buildings}.
Commentarii Math. Helv. {\bf 75 } (2000), 701-736

[CG] L. C{\sevenrm ARBONE}, H. G{\sevenrm ARLAND},
{\it Existence and deformations of lattices in groups of Kac-Moody type}.
In preparation

[Ch] F. C{\sevenrm HOUCROUN},
{\it Sous-groupes discrets des groupes $p$-adiques de rang un et 
arbres de Bruhat-Tits}.
Isr. J. Math. {\bf 93} (1996) 195-219

[Cho] A. C{\sevenrm HOSSON},
Isom\'etries dans les immeubles jumel\'es et construction d'automorphismes de
groupes de Kac-Moody.
PhD Univ. Amiens, 2000

[D] M. D{\sevenrm AVIS}, {\it Buildings are {\rm CAT(0)}}.
In Geometry and Cohomology in Group Theory, P.H. Kropholler, G.A. 
Niblo, R. St\"ohr
eds, LMS Lecture Notes Series {\bf 252} (1997) 108-123

[GP] D. G{\sevenrm ABORIAU}, F. P{\sevenrm AULIN},
{\it Sur les immeubles hyperboliques}. To appear in Geometri\ae $\,$ Dedicata

[H] J.-Y. H{\sevenrm \'EE}, {\it Construction de groupes tordus en
th\'eorie de Kac-Moody},
C. R. Acad. Sc. Paris {\bf 310} (1990) 77-80

[KP] V. K{\sevenrm AC}, D. P{\sevenrm ETERSON},
{\it Defining relations for certain infinite-dimensional groups}.
Ast\'erisque Hors-S\'erie (1984) 165-208

[Mar] G.A. M{\sevenrm ARGULIS}, Discrete subgroups of semisimple Lie groups.
Springer, 1990

[Mas] B. M{\sevenrm ASKIT}, Kleinian groups.
Grund. der Math. Wiss. {\bf 287}, Springer, 1987

[Mat] O. M{\sevenrm ATHIEU}, {\it Construction d'un groupe de 
Kac-Moody et applications}.
Compositio Math. {\bf 69} (1989) 37-60

[Mou] G. M{\sevenrm OUSSONG}, Hyperbolic Coxeter groups.
PhD Ohio State Univ., 1988

[MR] B. M{\sevenrm \"UHLHERR}, M. R{\sevenrm ONAN},
{\it Local to global structure in twin buildings}.
Invent. Math. {\bf 122} (1995) 71-81

[M\"u] B. M{\sevenrm \"UHLHERR},
{\it Locally split and locally finite buildings of $2$-spherical type}.
J. Reine angew. Math. {\bf 511} (1999) 119-143

[R\'e1] B. R{\sevenrm \'EMY},
Groupes de Kac-Moody d\'eploy\'es et presque d\'eploy\'es.
To appear in Ast\'erisque, Soci\'et\'e Math\'ematique de France, 2002

[R\'e2] B. R{\sevenrm \'EMY},
{\it Construction de r\'eseaux en th\'eorie de Kac-Moody}.
C. R. Acad. Sc. Paris {\bf 329} (1999) 475-478

[R\'e3] B. R{\sevenrm \'EMY},
{\it Classical and non-linearity properties of Kac-Moody lattices}.
To appear in Proc. \og Rigidity in Dynamics and
Geometry\fg (Newton Institute 2000), M. Burger and A. Iozzi eds, Springer, 2002

[R\'e4] B. R{\sevenrm \'EMY},
{\it Kac-Moody groups: split and relative theories. Lattices}.
To appear in Proc. \og Groups: Geometric and Combinatorial Aspects\fg
(Bielefeld 1999), London Math. Soc. Lecture Notes Series

[Ro1] M. R{\sevenrm ONAN},
Lectures on buildings.
Academic press, 1989

[Ro2] M. R{\sevenrm ONAN},
{\it Local isometries of twin buildings}.
Math. Zeitschrift {\bf 234} (2000) 435-455

[RT] M. R{\sevenrm ONAN}, J. T{\sevenrm ITS},
{\it Building buildings}.
Math. Ann. {\bf 278} (1987) 291-306

[S] J.-P. S{\sevenrm ERRE},
Arbres, amalgames, ${\rm SL}_2$.
Ast\'erisque {\bf 46}, 1977

[T1] J. T{\sevenrm ITS},
{\it Ensembles ordonn\'es, immeubles et sommes amalgam\'ees}.
Bull. Soc. Math. Belg. {\bf 38} (1986) 367-387

[T2] J. T{\sevenrm ITS},
{\it Uniqueness and presentation of Kac-Moody groups over fields}.
Journal of Algebra {\bf 105} (1987) 542-573

[T3] J. T{\sevenrm ITS}, {\it Th\'eorie des groupes}.
R\'esum\'e de cours, annuaire Coll\`ege de France (1989), 81-95

[T4] J. T{\sevenrm ITS}, {\it Th\'eorie des groupes}.
R\'esum\'e de cours, annuaire Coll\`ege de France (1990), 87-103

[T5] J. T{\sevenrm ITS}, {\it Twin buildings and groups of Kac-Moody type}.
In Groups, Combinatorics \& Geometry (Proc. LMS Symposium on Groups and
Combinatorics, Durham, July 1990), M. Liebeck and J. Saxl eds, LMS 
Lecture Notes Series
{\bf 165} (1992), Cambridge UP, 249-286

\vskip 10mm

I{\sevenrm NSTITUT } F{\sevenrm OURIER} -- UMR 5582
\hfill
D{\sevenrm EPARTMENT OF} M{\sevenrm ATHEMATICS},

U{\sevenrm NIVERSIT\'E} G{\sevenrm RENOBLE} 1
\hfill
S{\sevenrm TATISTICS AND} C{\sevenrm OMPUTER} S{\sevenrm CIENCE}

100, rue des maths -- BP 74
\hfill
U{\sevenrm NIVERSITY OF} I{\sevenrm LLINOIS AT} C{\sevenrm HICAGO}

38402 Saint-Martin-d'H\`eres, France
\hfill
851 S. Morgan Street

http:$/\!\!/$www-fourier.ujf-grenoble.fr/$^\sim$bremy
\hfill
Chicago, IL 60607-7045, USA

\pec

E-mail: {\tt bertrand.remy@ujf-grenoble.fr}
\hfill
E-mail: {\tt ronan@math.uic.edu}

\end

Bourse du British Council: CIBY D-TRA 6J1

---

322 Science \& Engineering Offices (SEO) m/c 249

---

[HKS] H{\sevenrm ERING}, K{\sevenrm ANTOR}, S{\sevenrm EITZ},
{\it Finite groups with a split $BN$-pair of rank $1$}. J. of Algebra {\bf 20}
(1972) 435-495

---

Auparavant en 4.E.2:

2) Once we have used the theorem of Gaboriau and Paulin, an 
alternative proof of 4.E.2 consists in showing
that the building $\Delta$ enjoys the Moufang property.
This property makes sense for single buildings -- see [Ro1, 6.4]; it 
involves automorphism groups indexed by
the root system $\Phi$ of the Weyl group.
What is left to do then is to use a general theorem due to V. 
Vermeulen [V], which says that a Moufang
building always comes from a group with a twin root datum.

---

At last, the additive parameters of the conjugating elements are
$\displaystyle {1 \over R} + {\mu^2 R^{-2} \over S -
\mu^2R^{-1}}$ and
$\displaystyle {1 \over R - \mu^2 S^{-1}}$.
Computing $\displaystyle {1 \over R - \mu^2 S^{-1}} - {1 \over
R}$, we get:

\pec

\centerline{$\displaystyle {R - R + \mu^2 S^{-1}
\over R \bigl( R - \mu^2 S^{-1} \bigr)}
= {\mu^2 S^{-1} \over R \bigl( R - \mu^2 S^{-1}
\bigr)}
\times {S \over S} = {\mu^2 \over R \bigl( R S -
\mu^2 \bigr)}
= {\mu^2R^{-1} \over R S - \mu^2}
= {\mu^2R^{-2} \over S - \mu^2R^{-1}}$,}

\pec

which proves the equality.

\pec

-----

Now we must compare two parameters, one for $u_a$ and one for $u_i$.
For $u_a$, after cancelling $k$, the exponent $-2\epsilon_a$, and 
dividing by $\lambda\mu$, this amounts to
checking that
$\displaystyle (R-{\mu^2\over S})S = R(S-{\mu^2\over R})$, which is obvious.
For $u_i$, after substracting $r$ and dividing by $-\lambda^2$, this 
amounts to showing that
$\displaystyle {1 \over R-{\mu^2S^{-1}}}={1 \over R}+{\mu^2 \over 
R^2(S-\mu^2R^{-1})}$,
but elementary algebra shows both sides equal
$\displaystyle {S \over RS-\mu^2}$.

--

This explains why we have to distinguish many cases.
For each of the five cases in the product formula, we must make specific
computations, which will be
referred to as {\it subcases}, one for each kind of generator $v$ of
$U^i$.
We have

\pec

-- the first subcase, when
$v=u_i(t)u_a(k)u_i(t)^{-1}$,
with $a \! \in \! {\cal P}(a_{1-i},i)$ and $t,k \! \in \! {\bf K}_i$;

-- the second subcase, when $v=u_i(t)u_a(k)u_i(t)^{-1}$,
with $a \! \in \! {\cal P}(a_{1-i},1-i)$, $t \! \in \! {\bf K}_i$ and 
$k \! \in \! {\bf K}_{1-i}$;

-- the third subcase, when $v=u_a(k)$,
with $a \! \in \! {\cal P}(a_i,i) \setminus \{ a_i \}$ and $k \! \in 
\! {\bf K}_i$;

-- the fourth subcase, when $v=u_a(k)$,
with $a \! \in \! {\cal P}(a_i,1-i)$ and $k \! \in \! {\bf K}_{1-i}$.

\pec

The list goes from the most complicated to the simplest case.
Indeed, the first two cases involve conjugating elements; in the 
first one, besides the action on the
indices of the root groups, both additive parameters are acted upon 
non-trivially.

--

This subsection explains why the construction of twin root data with 
several ground-fields
(hence the complete non-linearity phenomenon of Subsect. 5.B) is 
possible when the associated
buildings are right-angled Fuchsian.
In other words, we illustrate the following paradox: though they are 
two-dimensional, have Weyl groups of
arbitrarily large rank and share rigidity properties close to that of 
higher-rank symmetric spaces [BP], the
only local combinatorial structure of the buildings 
$I_{r,1+\underline q}$ makes them behave much like
trees, too.
The corollary of the proposition below shows that the flexibility of 
the links makes possible the
existence of characteristic 0 type lattices in a topological 
Kac-Moody group, which is itself close to an algebraic group over a
local field of positive characteristic (1.C.2).

-- 

If we specialize our study of Kac-Moody groups to the case of finite 
ground fields, we may
obtain groups such as ${\rm SL}_n({\bf F}_q[t,t^{-1}])$.
More generally, groups of points over Laurent polynomials of simply 
connected Chevalley groups are Kac-Moody (of
affine type), and at the same time are arithmetic groups over function fields.
This remark is the first step in seeing Kac-Moody groups over finite 
fields as discrete
groups.
More serious arguments for the analogy are given for instance in [A, 
CG, DJ, R\'e2, R\'e3].
In this paper, we are interested in the following further analogy.
Start again with the group $\Lambda:= {\rm SL}_n({\bf 
F}_q[t,t^{-1}])$: in the Kac-Moody context,
what is relevant is the diagonal $\Lambda$-action on the product of 
the Euclidean buildings
$\Delta_+$ and $\Delta_-$ attached to ${\rm SL}_n \bigl( {\bf F}_q 
(\!( t )\!) \bigr)$ and
${\rm SL}_n \bigl( {\bf F}_q (\!( t^{-1} )\!) \bigr)$, respectively.

--

Now, if we start with any abstract group $\Lambda$ enjoying the 
combinatorial structure of a twin
root datum, this group naturally acts diagonally on a pair of 
isomorphic twinned buildings.
Note that the group $\Lambda$ need not be obtained as a Kac-Moody group.
In the case where the buildings are locally finite, the closure $G$ 
of $\Lambda$ in
the full automorphism group of a single building is an interesting group.
Indeed, it acts (strongly transitively) on a building, it often 
contains natural lattices and,
mostly, is analogous to a semisimple group over a local field, as the 
above example of
${\rm SL}_n \bigl( {\bf F}_q (\!( t )\!) \bigr)$ shows.

--

In both theories the Weyl groups can be infinite and therefore do not 
provide a natural opposition relation, as given by
the longest word in a finite Weyl group.
Buildings with finite Weyl groups are called spherical, and the 
opposition relation is
crucial for Tits' classification of spherical buildings.
The theory of twin buildings, due initially to J. Tits and the second 
author, generalizes
the idea of a spherical building, by using a pair of buildings and 
allowing chambers in one building to be opposite chambers
in another building.

-- 

As a special case, consider the group $\Lambda:= {\rm SL}_n({\bf 
F}_q[t,t^{-1}])$.
It is a Kac-Moody group of affine type over the field ${\bf F}_q$ and 
is also a discrete arithmetic subgroup of
${\rm SL}_n \bigl( {\bf F}_q(t)\bigr)$.
It acts diagonally on the product of the affine buildings
$\Delta_+$ and $\Delta_-$ attached to ${\rm SL}_n \bigl( {\bf F}_q 
(\!( t )\!) \bigr)$ and
${\rm SL}_n \bigl( {\bf F}_q (\!( t^{-1} )\!) \bigr)$, respectively.
The buildings are locally finite because the base field ${\bf F}_q$ is finite.
If we forget $\Delta_-$, we get a non-discrete action of $\Lambda$ on 
$\Delta_+$ whose closure $G$ is
a topological group obtained as points of an algebraic group over a 
non-Archimedean local field,
namely ${\rm SL}_n \bigl( {\bf F}_q (\!( t )\!) \bigr)$.
On the other hand if we consider the stabilizer in $\Lambda$ of a 
vertex of $\Delta_-$, then we obtain the
subgroup ${\rm SL}_n({\bf F}_q[t^{-1}])$, which acts discretely on 
the other building $\Delta_+$.
More generally, one can start with any group $\Lambda$ admitting a 
twin root datum.
This could be a Kac-Moody group of some other type, not necessarily 
affine, but it could also be non-Kac-Moody, as the
present paper makes clear.
When the buildings $\Delta_+$ and $\Delta_-$ are locally finite, the 
closure $G$ of $\Lambda$ in
${\rm Aut}(\Delta_+)$ is analogous to the closure ${\rm SL}_n \bigl( 
{\bf F}_q (\!( t )\!) \bigr)$ of
${\rm SL}_n({\bf F}_q[t,t^{-1}])$.
We call such a group $G$ a {\it topological group of Kac-Moody type}, and more
specifically a {\it topological Kac-Moody group} if $\Lambda$ is 
actually a Kac-Moody group over a
finite field.

--

For any $n,n' \geq 0$, the pair $\{a_i(n);a_j(n')\}$ of positive 
roots fails to be prenilpotent because
the intersection of the negative roots $-a_i$ and $-a_j$, hence also 
of $-a_i(n)$ and $-a_j(n')$, is
empty.

--

The Coxeter complex of $(W,S)$ admits a simplicial realization where 
the maximal simplices,
representing the chambers, are acted upon simply transitively by $W$ 
[Ro1 \S 2].
A root $a$ is represented by a half space whose boundary $\partial a$ 
is called a {\it wall} [loc. cit.].
Walls are thus in bijection with pairs of opposite roots.